\pdfoutput=1 
\documentclass{article}
\usepackage{excludeonly}

\usepackage{savesym}
\usepackage{amsthm,amsmath}
\usepackage{amssymb,latexsym,graphicx}
\usepackage{amscd}
 \usepackage[all,cmtip]{xy}
\usepackage{tikz-cd}
\usepackage{accents}
\usepackage{cite}
\usepackage{mathtools}

\usepackage{hycolor}
\usepackage{xcolor}
\usepackage[
                       colorlinks=true,
                       linkcolor=black, %red,
                       citecolor=black, %green,
                       urlcolor=blue,
                     ]{hyperref}
\usepackage{booktabs} % tables

\newtheorem{theoremABC}{Theorem}

\newtheorem{theorem}{Theorem}[section]

\newtheorem{corollary}[theorem]{Corollary}

\newtheorem{lemma}[theorem]{Lemma}
\newtheorem{proposition}[theorem]{Proposition}

\theoremstyle{definition}
\newtheorem{definition}[theorem]{Definition}

\newtheorem{remark}[theorem]{Remark}

\theoremstyle{remark}
\newtheorem*{notation}{Notation}
\newcommand{\C}{{\mathbb{C}}}

\newcommand{\LL}{{\mathbb{L}}}

\newcommand{\N}{{\mathbb{N}}}

\newcommand{\Q}{{\mathbb{Q}}}
\newcommand{\R}{{\mathbb{R}}}
\renewcommand{\SS}{{\mathbb{S}}}

\newcommand{\Z}{{\mathbb{Z}}}
\newcommand{\Aa}{{\mathcal{A}}}   % connections
\newcommand{\Bb}{{\mathcal{B}}}
\newcommand{\Hh}{{\mathcal{H}}}

\newcommand{\Ll}{{\mathcal{L}}}   % Lagrangian planes
\newcommand{\Qq}{{\mathcal{Q}}}
\newcommand{\Rr}{{\mathcal{R}}}
\newcommand{\Ss}{{\mathcal{S}}}
\newcommand{\Tt}{{\mathcal{T}}}
\newcommand{\Uu}{{\mathcal{U}}}

\newcommand{\dom}{{\rm dom\, }}           % domain
\newcommand{\grad}{\mathop{\mathrm{grad}}}    % gradient
\newcommand{\CZ}{{\mu_{\rm CZ}}}              % Conley-Zehnder index
\newcommand{\CZcan}{{\mu^{\rm CZ}}}           % clockw. Conley-Zehnder index
\newcommand{\cgraph}[1]{\Gamma_{\kern-.5ex{}#1}}     % contact graph map
\newcommand{\ind}{{\rm Ind}} 
\newcommand{\Crit}{{\rm Crit}}        % Critical points
\newcommand{\CAP}{\mathop{\cap}}           % cap in formulas: adjusted spaces
 
\newcommand{\spec}{\mathrm{spec}\,}        % spectrum
\newcommand{\norm}{{\rm norm}}

\newcommand{\inner}[2]{\langle #1, #2\rangle}   % inner product < , >
\newcommand{\INNER}[2]{\left\langle #1, #2\right\rangle}

\def\NABLA#1{{\mathop{\nabla\kern-.5ex\lower1ex\hbox{$#1$}}}}
\def\Nabla#1{\nabla\kern-.5ex{}_{#1}}
\def\Tabla#1{\Tilde\nabla\kern-.5ex{}_{#1}}
\def\abs#1{\mathopen|#1\mathclose|}   
\def\Abs#1{\left|#1\right|}            
\def\norm#1{\mathopen\|#1\mathclose\|}
\def\Norm#1{\left\|#1\right\|}

\renewcommand{\Tilde}{\widetilde}

\newcommand{\p}{{\partial}}

\newcommand{\INTO}{\hookrightarrow}              % embedding
\renewcommand{\1}{{{\mathchoice {\rm 1\mskip-4mu l} {\rm 1\mskip-4mu l}
{\rm 1\mskip-4.5mu l} {\rm 1\mskip-5mu l}}}}
\newlength\eqshift
\renewcommand\theequation{\thesection.\arabic{equation}}
\let\savetheequation\theequation

\makeatletter
\renewcommand*\env@matrix[1][\arraystretch]{%
  \edef\arraystretch{#1}%
  \hskip -\arraycolsep
  \let\@ifnextchar\new@ifnextchar
  \array{*\c@MaxMatrixCols c}}
\makeatother

\makeatletter
\let\save@mathaccent\mathaccent
\newcommand*\if@single[3]{%
  \setbox0\hbox{${\mathaccent"0362{#1}}^H$}%
  \setbox2\hbox{${\mathaccent"0362{\kern0pt#1}}^H$}%
  \ifdim\ht0=\ht2 #3\else #2\fi
  }
\newcommand*\rel@kern[1]{\kern#1\dimexpr\macc@kerna}
\newcommand*\widebar[1]{\@ifnextchar^{{\wide@bar{#1}{0}}}{\wide@bar{#1}{1}}}
\newcommand*\wide@bar[2]{\if@single{#1}{\wide@bar@{#1}{#2}{1}}{\wide@bar@{#1}{#2}{2}}}
\newcommand*\wide@bar@[3]{%
  \begingroup
  \def\mathaccent##1##2{%
    \let\mathaccent\save@mathaccent
    \if#32 \let\macc@nucleus\first@char \fi
    \setbox\z@\hbox{$\macc@style{\macc@nucleus}_{}$}%
    \setbox\tw@\hbox{$\macc@style{\macc@nucleus}{}_{}$}%
    \dimen@\wd\tw@
    \advance\dimen@-\wd\z@
    \divide\dimen@ 3
    \@tempdima\wd\tw@
    \advance\@tempdima-\scriptspace
    \divide\@tempdima 10
    \advance\dimen@-\@tempdima
    \ifdim\dimen@>\z@ \dimen@0pt\fi
    \rel@kern{0.6}\kern-\dimen@
    \if#31
      \overline{\rel@kern{-0.6}\kern\dimen@\macc@nucleus\rel@kern{0.4}\kern\dimen@}%
      \advance\dimen@0.4\dimexpr\macc@kerna
      \let\final@kern#2%
      \ifdim\dimen@<\z@ \let\final@kern1\fi
      \if\final@kern1 \kern-\dimen@\fi
    \else
      \overline{\rel@kern{-0.6}\kern\dimen@#1}%
    \fi
  }%
  \macc@depth\@ne
  \let\math@bgroup\@empty \let\math@egroup\macc@set@skewchar
  \mathsurround\z@ \frozen@everymath{\mathgroup\macc@group\relax}%
  \macc@set@skewchar\relax
  \let\mathaccentV\macc@nested@a
  \if#31
    \macc@nested@a\relax111{#1}%
  \else
    \def\gobble@till@marker##1\endmarker{}%
    \futurelet\first@char\gobble@till@marker#1\endmarker
    \ifcat\noexpand\first@char A\else
      \def\first@char{}%
    \fi
    \macc@nested@a\relax111{\first@char}%
  \fi
  \endgroup
}
\makeatother

\long\def\symbolfootnote[#1]#2{\begingroup%
\def\thefootnote{\fnsymbol{footnote}}\footnote[#1]{#2}\endgroup}
\begin{document}
\sloppy

\author{\quad Urs Frauenfelder \quad \qquad\qquad
             Joa Weber\footnote{
  Email: urs.frauenfelder@math.uni-augsburg.de
  \hfill
  joa@ime.unicamp.br
  }
    \\
    Universit\"at Augsburg \qquad\qquad%\quad
    UNICAMP
}

\title{The regularized free fall \\ I -- Index computations}
\date{\today}

\maketitle

%%%%%%%%%%%%%%%%%%%%%%%%%%%%%%%%%%%
%%%%%%% Abstract %%%%%%%%%%%%%%%%%%%%%
%%%%%%%%%%%%%%%%%%%%%%%%%%%%%%%%%%%
\begin{abstract}
Main results are, firstly, a generalization of the Conley-Zehnder index from
ODEs to the delay equation at hand and, secondly, the equality of the Morse index
and the clockwise normalized Conley-Zehnder index $\CZcan$.

We consider the non-local Lagrangian action functional $\Bb$
discovered by Barutello, Ortega, and Verzini~\cite{Barutello:2021b}
with which they obtained a new regularization of the Kepler problem.
Critical points of this functional are regularized periodic
solutions $x$ 
of the Kepler problem. In this article
we look at \emph{period 1 only} and at dimension one (gravitational free fall).

Via a non-local Legendre transform regularized periodic Kepler orbits $x$
can be interpreted as periodic solutions $(x,y)$ of a Hamiltonian delay equation.
In particular, regularized $1$-periodic solutions of the free fall are
represented variationally in two ways: As critical points $x$ of
a non-local Lagrangian action functional and as critical points $(x,y)$
of a non-local Hamiltonian action functional. 
\\
As critical points of the Lagrangian action
the $1$-periodic solutions have a finite Morse index which we compute first.
\\
As critical points of the Hamiltonian action $\Aa_\Hh$ one encounters
the obstacle, due to non-locality, that the $1$-periodic
solutions are not generated any more by a flow on the
phase space manifold. Hence the usual definition of the Conley-Zehnder
index as intersection number with a Maslov cycle is not available.
In the local case Hofer, Wysocki, and Zehnder~\cite{Hofer:1995b} gave an
alternative definition of the Conley-Zehnder index by assigning a winding
number to each eigenvalue of the Hessian of $\Aa_\Hh$ at critical points.

In this article we show how to generalize the Conley-Zehnder index to
the non-local case at hand.
On one side we discover how properties from the local case
generalize to this delay equation, and on the other side
we see a new phenomenon arising.
In contrast to the local case the winding number is not any more monotone
as a function of the eigenvalues.
\end{abstract}

\newpage
\tableofcontents

\newpage
%%%%%%%%%%%%%%%%%%%%%%%%%%%%%%%%%%%
%%%%%%%%%%%%%%%%%%%%%%%%%%%%%%%%%%%
%%%%%%% Section:  %%%%%%%%%%%%%%%%%%%%%
%%%%%%%%%%%%%%%%%%%%%%%%%%%%%%%%%%%
%%%%%%%%%%%%%%%%%%%%%%%%%%%%%%%%%%%
\section{Introduction}

In celestial mechanics, as well as in the theory of atoms,
collisions play an intriguing role.
There are many geometric ways how to regularize collisions,
see for instance~\cite{Levi-Civita:1920a,Moser:1970a}.
In both regularizations, Levi-Civita and Moser, one has to change
time.
A quite new approach to the regularization of collisions was
discovered in the recent article~\cite{Barutello:2021b}
by Barutello, Ortega, and Verzini where the change of time
leads to a \emph{\color{brown} delayed} functional
(meaning that the critical point equation is a delay equation).
In Section~\ref{sec:ff} we explain, starting from the physics
1-dimensional Kepler problem, how one arrives at this functional
\begin{equation*}
\begin{split}
   \Bb\colon W^{1,2}_\times:= W^{1,2}(\SS^1,\R)\setminus \{0\}&\to(0,\infty) \\
   x&\mapsto 
   {\color{brown} 4 \Norm{x}^2} \tfrac12 \norm{ x^\prime}^2+\frac{1}{\norm{x}^2}
\end{split}
\end{equation*}
where $\norm{\cdot}$ is the $L^2$ norm associated to the $L^2$ inner
product $\inner{\cdot}{\cdot}$.
One might interpret the functional $\Bb$
as a \textbf{non-local Lagrangian action functional},
a non-local mechanical system, consisting of
kinetic minus potential energy
\[
   \Bb(x)=\tfrac12\inner{ x^\prime}{ x^\prime}_x+\frac{1}{\norm{x}^2}
\]
Here we use the following metric on the tangent bundle of the loop space
\[
   \inner{\cdot}{\cdot}_x:= {\color{brown} 4 \Norm{x}^2}
   \inner{\cdot}{\cdot}
\]
which from the perspective of the manifold is non-local as it depends on
the whole loop. Also the potential $x\mapsto -\frac{1}{\norm{x}^2}$
is only defined on the loop space.
\\
The set $\Crit\,\Bb$ of critical points of the non-local action
consists of the smooth solutions
$x\in C^\infty(\SS^1,\R)\setminus \{0\}$ of the second order delay,
or non-local, equation
\begin{equation}\label{eq:Crit-Lag}
   x^{\prime\prime}
   =\alpha x,\qquad
   \alpha=\alpha_x :=\left(
   \frac{\norm{ x^\prime}^2}{\norm{x}^2}-\frac{1}{2\norm{x}^6}
   \right) <0
\end{equation}

Nevertheless, this non-local Lagrangian action functional admits a
Legendre transform which leads to a non-local Hamiltonian that lives
on (the $L^2$ extension of) the cotangent bundle of the loop space, namely
\begin{equation*}
\begin{split}
   \Hh\colon W^{1,2}_\times \times L^{2}&\to\R \\
   (x,y)&\mapsto
   \tfrac12\inner{y}{y}^x -\frac{1}{\norm{x}^2}
\end{split}
\end{equation*}
where we use the dual metric on the cotangent bundle of the loop space
\[
   \inner{\cdot}{\cdot}^x:= {\color{brown} \tfrac{1}{4 \Norm{x}^2}}
   \inner{\cdot}{\cdot}
\]
Associated to the non-local Hamiltonian $\Hh$ there is the
\textbf{non-local Hamiltonian action functional}, namely
\begin{equation*}
\begin{split}
   \Aa_\Hh:=\Aa_0-\Hh\colon 
   W^{1,2}_\times \times L^{2}&\to\R\\
   (x,y)&\mapsto\int_0^1y(\tau) x^\prime(\tau)\, d \tau-\Hh(x,y)
\end{split}
\end{equation*}
The natural base point projection and the following injection
\begin{equation}\label{eq:proj-inj}
\begin{aligned}
   \pi\colon W^{1,2}_\times\times L^{2}&\to W^{1,2}_\times\qquad&
   \iota\colon W^{1,2}_\times&\to W^{1,2}_\times\times L^{2}
   \\
   (x,y)&\mapsto x\qquad&
   x&\mapsto\left( x,{\color{brown} 4\norm{x}^2}x^\prime\right)
\end{aligned}
\end{equation}
are \emph{along the sets of critical points} inverses of one another
-- which of course explains the choice of the factor
${\color{brown} 4\norm{x}^2}$. Therefore the sets
\begin{equation*}
\begin{tikzcd} [row sep=tiny] %[column sep=small] %
\Crit\,\Aa_\Hh
\arrow[rr, shift left=2,"\pi", "1:1"']
  &&\Crit\,\Bb
      \arrow[ll, shift left=2,"\iota"]
\end{tikzcd} 
\end{equation*}
are in one-to-one correspondence.
Whereas $\Bb$ at any critical point has a finite dimensional Morse index,
in sharp contrast both the Morse index and coindex are infinite at the
critical points of $\Aa_\Hh$.

Since the action functional $\Aa_\Hh$ is non-local its critical points
cannot be interpreted as fixed points of a flow.
Therefore the usual definition of the Conley-Zehnder index
as a Maslov index does not make sense.
However, in~\cite{Hofer:1995b} Hofer, Wysocki, and Zehnder
gave a different characterization of the Conley-Zehnder index
in terms of winding numbers of the eigenvalues of the Hessian.

We show in this article that this definition
of the Conley-Zehnder index makes sense, too,
for the critical points of the \emph{non-local} functional $\Aa_\Hh$.
Our main result is equality of Morse and Conley-Zehnder index
of corresponding critical points.

\begin{theoremABC}
For each critical point $(x,y)$ of $\Aa_\Hh$ the 
canonical (clockwise normalized) Conley-Zehnder index
\[
    \CZcan(x,y)=\ind(x)
\]
is equal to the Morse index of $\Bb$ at $x$.
\end{theoremABC}

\begin{proof}
Proposition~\ref{prop:Morse-calc} and
Proposition~\ref{prop:CZ-calc}.
\end{proof}

The theorem generalizes the local result, see~\cite{weber:2002a}, to this
\emph{non-local} case. Note that~\cite[Thm.~1.2]{weber:2002a} uses the
counter-clockwise normalized Conley-Zehnder index $\CZ=-\CZcan$.
Sign conventions are discussed at large in the
introduction to~\cite{Weber:2017b}.
The local result was a crucial ingredient
in the proof that Floer homology of the cotangent bundle
is the homology of the loop
space~\cite{Viterbo:1998a,salamon:2006a,abbondandolo:2006b}.

\smallskip
\begin{notation}
We define $\SS^1:=\R/\Z$ and consider functions $x\colon\SS^1\to\R$
as functions defined on $\R$ that satisfy $x(\tau+1)=x(\tau)$ for
every $\tau\in\R$.
Throughout $\inner{\cdot}{\cdot}$ is the standard $L^2$ inner product
on $L^2(\SS^1,\R)$ and $\norm{\cdot}$ is the induced $L^2$ norm.
\end{notation}

\smallskip\noindent
{\bf Outlook.}
Since 2018 first steps have been taken to study Floer homology for delay
equations, see~\cite{Albers:2020a,Albers:2018c,Albers:2018a}.
The present article is the first in a series of four dealing with
the free fall -- the simplest instance
which might already exhibit all novelties that occur in
comparison to ode Floer homology of the cotangent bundle
(which still represents the homology of a space, loop space).
The other parts will deal with II ``Homology computation via heat flow'',
III ``Floer homology'', and IV ``Floer homology and heat flow homology''.

\smallskip\noindent
{\bf Acknowledgements.}
UF acknowledges support by DFG grant
FR 2637/2-2.

%%%%%%%%%%%%%%%%%%%%%%%%%%%%%%%%%%%
%%%%%%%%%%%%%%%%%%%%%%%%%%%%%%%%%%%
%%%%%%% Section:  %%%%%%%%%%%%%%%%%%%%%
%%%%%%%%%%%%%%%%%%%%%%%%%%%%%%%%%%%
%%%%%%%%%%%%%%%%%%%%%%%%%%%%%%%%%%%
\section{Free gravitational fall}\label{sec:ff}

Section~\ref{sec:ff} is to motivate and explain regularization. Readers familiar
with regularization could go directly to subsequent sections.

%%%%%%%%%%%%%%%%%%%%%%%%%%%%%%%%%%%
%%%%%%% Subsection:  %%%%%%%%%%%%%%%%%%
%%%%%%%%%%%%%%%%%%%%%%%%%%%%%%%%%%%
\subsection{Classical}

For $r>0$ and $\frak v\in\R$ let $L(r, \frak v):=\frac12\abs{\frak v}^2-V(r)$ where
$V(r):=-1/r$. Then
\[
   d_{\frak v} L(r, \frak v)=\frak v ,\qquad d_rL(r, \frak v)=-1/r^2
\]
Due to the potential one cannot allow $r=0$.
Thus the classical action functional
\begin{equation}\label{eq:class-action}
   \Ss_L(q):=\int_0^1 L(q(t),\dot q(t))\, dt
   =\int_0^1\left(\tfrac12\Abs{\dot q(t)}^2+\frac{1}{q(t)}\right) dt
\end{equation}
is defined on the
space
\begin{equation}\label{eq:W12+}
   W^{1,2}_{+}:=W^{1,2}(\SS^1,(0,\infty))
\end{equation}
that consists of absolutely continuous positive functions
$q\colon[0,1]\to(0,\infty)$ which are periodic, that is $q(1)=q(0)$,
and whose derivative is $L^2$ integrable, that is
$\norm{\dot q}^2:=\int_0^1 \dot q(t)^2\, dt<\infty$.
\\
The Euler-Lagrange (or critical point) equation is given by the second order ode
\begin{equation}\label{eq:EL}
   \frac{d}{dt} d_2 L(q,\dot q)=d_1 L(q,\dot q)
   \qquad\Leftrightarrow\qquad
   \ddot q=-\frac{1}{q^2}
\end{equation}
pointwise at $t$.
Unfortunately, there are no periodic
solutions of this equation. 
In other words, there are no critical points of $\Ss_L$ on the space
$W^{1,2}_+$, in symbols
\begin{equation}\label{eq:Crit-Ss-empty}
   \Crit\,\Ss_L=\emptyset
\end{equation}
All solutions $q$ with zero initial velocity $v_0$
end up in collision with the origin.
Set $q_0:=q(0)$ and $v_0:=\dot q(0)$. Multiply~(\ref{eq:EL}) by
$2\dot q$ and integrate to get
$$
   \dot q(t)^2-v_0^2
   =\int_0^t \underbrace{\ddot q 2\dot q}
                     _{\frac{d}{ds} {\dot q}^2} ds
   \stackrel{(\ref{eq:EL})}{=}\int_0^t \underbrace{- \frac{2\dot q}{q^2}}
                        _{2\frac{d}{ds} {\dot q}^{-1}} ds
   =\frac{2}{q(t)}-\frac{2}{q_0}
$$
hence
$$
   \dot q(t)=\pm\sqrt{\frac{2}{q(t)}+2E},\qquad
   E:=\frac{v_0^2}{2}-\frac{1}{q_0}\equiv \frac{\dot q(t)^2}{2}-\frac{1}{q(t)}
$$
This shows that collision and bouncing off happen at infinite velocity.
If the particle starts with zero initial velocity,
collision happens in finite time since the acceleration
$\ddot q$ is bounded away from zero.

%\newpage
%%%%%%%%%%%%%%%%%%%%%%%%%%%%%%%%%%%
%%%%%%% Subsection:  %%%%%%%%%%%%%%%%%%
%%%%%%%%%%%%%%%%%%%%%%%%%%%%%%%%%%%
\subsection{Non-local regularization}

In this section we explain why Barutello, Ortega, and
Verzini~\cite{Barutello:2021b} discovered their functional $\Bb$
defined on the larger (than $W^{1,2}_+$) space
\[
   W^{1,2}_\times:=W^{1,2}_\times(\SS^1,\R)\setminus \{0\}
   \quad\supset\quad
   W^{1,2}(\SS^1,(0,\infty))=:W^{1,2}_+
\]
It consists of absolutely continuous maps $x\colon[0,1]\to\R$
which are periodic $x(1)=x(0)$, not identically zero $x\not\equiv 0$,
and whose weak derivative $ x^\prime$ is $L^2$ integrable. 

The key observation is that if on the subset $W^{1,2}_+$ -- the
domain~(\ref{eq:W12+}) of the classical action $\Ss_L$ --
one defines an operation $\Rr$ that
takes the square and appropriately rescales time, then
the composition $(\Ss_L\circ \Rr)(x)$ is given by a formula,
let's call it $\Bb(x)$, see Figure~\ref{fig:fig-Q^*S},
that makes perfectly sense on the ambient space
$W^{1,2}_\times$ whose elements are real-valued and so allow
for origin traversing.

Whereas the classical functional $\Rr^*\Ss_L$ is defined on loops
that take values in $(0,\infty)$, thereby not allowing for collisions at $0$,
has no critical points, the rescaled functional has critical points
when considered on the extended domain
$W^{1,2}_\times(\SS^1,\R)\setminus \{0\}$.
These critical points correspond to classical solutions having
collisions, that is running into the origin $0$.

\begin{figure}[h]
  \centering
\begin{equation*}
\begin{tikzcd} [row sep=tiny] %[column sep=small] %
\R
  &\R
  &&\R
\\ \mbox{} \\ \mbox{} \\
W^{1,2}_\times 
\arrow[uuu, "\Bb"]
\arrow[r, phantom, "\supset"]
  &W^{1,2}_+  
    \arrow[rr, "\Rr", "1:1"']
    \arrow[uuu, dashed, "\Bb:=", "\Rr^*\Ss_L"']
  &&W^{1,2}_+ 
      \arrow[uuu, "\Ss_L"]
\\
  &x
    \arrow[rr, mapsto]
  &&
    q_x:=x^2\circ \tau_x
\\
  &{\color{gray} \text{time $\tau$}}
  &&
    {\color{gray} \text{time $t$}}
\end{tikzcd} 
\end{equation*}
  \caption{Pull-back $\Rr^*\Ss_L$ gives a formula $\Bb$ that makes
                 sense on a larger space on which critical points exist
                 and represent time rescaled physical trajectories}
   \label{fig:fig-Q^*S}
\end{figure}

The key fact is that $\Bb$ admits critical points $x$ on the ambient
space $W^{1,2}_\times$, see~(\ref{eq:q-heat}), in which case
$q_x:=\Rr(x)$ solves by Proposition~\ref{prop:reg=phys}
the classical free fall equation~(\ref{eq:EL})
away from a finite set $T_x\subset\SS^1$ of collision times.

\boldmath
%%%%%%% Subsubsection:  %%%%%%%%%%%%%%%%
\subsubsection*{Time rescaling on the small space $W^{1,2}_+$}
\unboldmath

Pick a loop $x\in W^{1,2}_+$, see~(\ref{eq:W12+}).
We call the variable of the map $x\colon[0,1]\to(0,\infty)$
the \textbf{regularized time}, usually denoted by $\tau$.
\textbf{Classical time} we call the values of the map
$t_x\colon[0,1]\to[0,1]$ defined by
\begin{equation}\label{eq:class-time-tau}
   t_x(\tau)
   :=\frac{\int_0^\tau x(s)^2\, ds}{\norm{x}^2}
\end{equation}
Classical time has the following properties
\begin{equation}\label{eq:class-time-tau-deriv}
   t_x^\prime(\tau)
   =\frac{x(\tau)^2}{\norm{x}^2}> 0
   ,\quad
   t_x\in C^1,\quad t_{r x}=t_x,\quad
   t_x(0)=0,\quad t_x(1)=1
\end{equation}
for every real $r>0$.
Since, moreover, the map $t_x\colon[0,1]\to[0,1]$ is strictly
increasing, it is a bijection and we denote its inverse by
\[
   \tau_x\colon[0,1]\to[0,1]
\]
The inverse $\tau_x:={t_x}^{-1}$ of classical time 
inherits the property $\tau_{rx}=\tau_x$ $\forall r>0$~and
\begin{equation}\label{eq:inverse}
   \dot \tau_x(t)
   =\frac{1}{t_x^\prime(\tau_x(t))}
   =\frac{\norm{x}^2}{x(\tau_x(t))^2}
   ,\quad
   \tau_x\in C^1,\quad
   \tau_x(0)=0,\quad \tau_x(1)=1
\end{equation}

\begin{definition}\label{def:reg-path}
The \textbf{rescale-square operation}
is defined by
\[
   \Rr\colon W^{1,2}_+\to W^{1,2}_+
   ,\quad
   x\mapsto x^2\circ \tau_x
\]
We abbreviate $q_x:=\Rr(x)$. 
Note that $\Rr(rx)=r^2\Rr(x)$ for $r>0$.
\end{definition}

\begin{lemma}[Well defined bijection]
For $x\in W^{1,2}_+$ the image $\Rr(x)$ lies in $W^{1,2}_+$, too,
and the map $\Rr\colon W^{1,2}_+\to W^{1,2}_+$ is bijective with
inverse~(\ref{eq:Rr-inverse}).
\end{lemma}

\begin{proof}
We must show that both are in $L^2$, namely
a) $q_x(t)=x^2(\tau_x(t))$ and b)
\[
   \dot q_x(t)
   =2x(\tau_x(t))\, x^\prime(\tau_x(t))\, \dot\tau_x(t)
   \stackrel{(\ref{eq:inverse})}{=}2\norm{x}^2\,\frac{x^\prime(\tau_x(t))}{x(\tau_x(t))}
\]
Use that $x$ is continuous,
so $\norm{x}_{L^\infty}<\infty$, to get a)
\[
   \norm{q_x}^2
   =\int_0^1 x^4(\tau_x(t))\, dt
   =\norm{x}_{L^\infty}^4<\infty
\]
In what follows change the variable to $\sigma:=\tau_x(t)$
using~(\ref{eq:inverse}) to get b)
\begin{equation}\label{eq:fsfcs}
   \norm{\dot q_x}^2
   =\int_0^1\left(2\norm{x}^2\,\frac{x^\prime(\tau_x(t))}{x(\tau_x(t))}\right)^2
         \frac{x^2(\sigma)\, d\sigma}{\norm{x}^2}
   =4\norm{x}^2\INNER{x^\prime}{x^\prime}<\infty
\end{equation}
Here the value is finite since $x$ is of class $W^{1,2}$.
This proves that $\Rr(x)\in W^{1,2}_+$.

\smallskip\noindent
{\sc Surjective.}
Given $q\in W^{1,2}_+$, set
\begin{equation}\label{eq:Rr-inverse}
   \Qq(q):=x_q:=q^\frac12\circ \tau_{1/\sqrt{q}}
\end{equation}
Then for $\tau\in[0,1]$ we obtain
\[
   t_{x_q}(\tau)
   \stackrel{\text{def.}}{=}
   \frac{\int_0^\tau{x_q}^2(\sigma)\, d\sigma}{\int_0^1{x_q}^2(\sigma)\, d\sigma}
   =\frac{\int_0^\tau q\circ\overbrace{\tau_{1/\sqrt{q}}(\sigma)}^{=:s}\, d\sigma}
      {\int_0^1 q\circ \tau_{1/\sqrt{q}}(\sigma)\, d\sigma}
   =\frac{\int_0^{\tau_{1/\sqrt{q}}(\tau)} \frac{ds}{\norm{1/\sqrt{q}}^2}}
      {\int_0^1 \frac{ds}{\norm{1/\sqrt{q}}^2}}
   =\tau_{1/\sqrt{q}}(\tau)
\]
by change of variables. With this result we get that
\[
   \bigl(\Rr\circ\underbrace{\Qq(q)}_{x_q}\bigr)(t)
%   =\left(\Rr(q_x)\right)(t)
   \stackrel{\text{def. $\Rr$}}{=}
   x_q^2\circ \tau_{x_q}(t)
   \stackrel{\text{def. $x_q^2$}}{=}
   q\circ\underbrace{t_{1/\sqrt{q}}}_{t_{x_q}}\circ
   \underbrace{\tau_{x_q}}_{{t_{x_q}}^{-1}}(t)
   =q(t)
\]
{\sc Injective.}
For $x\in W^{1,2}_+$ set
$q_x:=\Rr(x):=x^2\circ \tau_x$. Then for $t\in[0,1]$ we get
\[
   \int_0^t\frac{1}{q_x(s)}\,ds
   =\int_0^t\frac{1}{x^2\circ \tau_x(s)}\,ds
   =\int_0^{\tau_x(t)}\frac{1}{x^2(\sigma)}\frac{x^2(\sigma)\, d\sigma}{\norm{x}^2}
   =\frac{\tau _x(t)}{\norm{x}^2}
\]
by change of variables $\sigma=\tau_x(s)$ using~(\ref{eq:inverse}).
Pick $t=1$ to obtain
\begin{equation}\label{eq:1/q^2}
   \int_0^1\frac{1}{q_x(t)}\, dt
   =\frac{1}{\norm{x}^2}
\end{equation}
Therefore we get for $t\in[0,1]$ the formula
\[
   \tau_x(t)
   =\norm{x}^2 \int_0^t\frac{1}{q_x(s)}\,ds
   \stackrel{(\ref{eq:1/q^2})}{=}
      \frac{\int_0^t\frac{1}{q_x(s)}\,ds}
      {\int_0^1\frac{1}{q_x(s)}\,ds}
   \stackrel{(\ref{eq:class-time-tau})}{=:}
   t_{1/\sqrt{q_x}}(t)
\]
With this result we obtain
\[
   \bigl(\Qq\circ\underbrace{\Rr(x)}_{q_x}\bigr)(\tau)
   \stackrel{\text{def. $\Qq$}}{=}
      \sqrt{q_x}\circ \tau_{1/\sqrt{q_x}}(\tau)
   \stackrel{\text{def. $q_x$}}{=}
      x\circ\underbrace{\tau_x}_{t_{1/\sqrt{q_x}}}\circ
      \underbrace{t_{1/\sqrt{q_x}}}_{{\tau_{1/\sqrt{q_x}}}^{-1}}(\tau)
   =x(\tau)
\]
\end{proof}

\begin{lemma}[Pull-back yields an extendable formula]
Given $x\in W^{1,2}_+$, then
\[
   \Ss_L(\Rr(x))
   =\tfrac12\inner{x^\prime}{x^\prime}_x+\frac{1}{\norm{x}^2}
\]
as illustrated by Figure~\ref{fig:fig-Q^*S}.
\end{lemma}

Note that in the previous formula $x(\tau)$ may take on the value zero at will,
even along intervals, as long as $\norm{x}\not=0$, i.e. as long as $x$ is not
constantly zero.

\begin{proof}
Set $q_x:=\Rr(x)$.
Definition~(\ref{eq:class-action}) of $\Ss_L$ tells that
\begin{equation}\label{eq:Ss=Bb}
\begin{split}
   \Ss_L(q_x)
   \stackrel{(\ref{eq:class-action})}{=}
   \tfrac12\norm{\dot q_x}^2+
   \int_0^1\frac{1}{q_x(t)}\, dt
   =\tfrac12\, 4\norm{x}^2\INNER{x^\prime}{x^\prime}+\frac{1}{\norm{x}^2}
\end{split}
\end{equation}
where in the second equality we used~(\ref{eq:fsfcs}) and~(\ref{eq:1/q^2}).
\end{proof}

%\newpage
%%%%%%%%%%%%%%%%%%%%%%%%%%%%%%%%%%%
%%%%%%% Subsection:  %%%%%%%%%%%%%%%%%%
%%%%%%%%%%%%%%%%%%%%%%%%%%%%%%%%%%%
\subsection{Correspondence of solutions}

Whereas the classical action $\Ss_L$ does not admit, see~(\ref{eq:Crit-Ss-empty}),
on the small space $W^{1,2}_+$ any critical points, that is no
$1$-periodic free fall solutions, the formula~(\ref{eq:Ss=Bb}) for
$\Ss_L$ evaluated on the regularization $q_x:=\Rr(x)$ of a loop
$x\in W^{1,2}_+$ makes perfectly sense on the larger Sobolev Hilbert space
$W^{1,2}_\times:=W^{1,2}(\SS^1,\R)\setminus\{0\}$ with the origin removed. 
This way one arrives, in case of the free fall, at the
Barutello-Ortega-Verzini~\cite{Barutello:2021b} non-local functional
\begin{equation}\label{eq:BOV-fctl}
\begin{split}
   \Bb\colon W^{1,2}_\times:=W^{1,2}(\SS^1,\R)\setminus \{0\}\to(0,\infty)
   ,\quad
   x\mapsto \tfrac12\INNER{x^\prime}{x^\prime}_x+\frac{1}{\norm{x}^2}
\end{split}
\end{equation}
For this functional $\Bb$ eventual zeroes of $x$ cause no problem at all.
It has even lots of critical points - one circle worth of critical
points for each natural number $k\in\N$, see Lemma~\ref{le:crit-pts-Bb}.

But are the critical points $x\colon\SS^1\to\R$ of $\Bb$
related to free fall solutions $q\colon \dom q\to(0,\infty)$?
If so what is the domain of~$q$?
Next we prepare to give the answers in Proposition~\ref{prop:reg=phys}
below. From now on
\begin{itemize}\setlength\itemsep{0ex} 
\item
  we only consider critical points $x\in W^{1,2}_\times$ of $\Bb$ and
\item
  we identify $\SS^1$ with $[0,1]/\{0,1\}$. 
\end{itemize}
A priori such $x\colon[0,1]\to\R$
might have zeroes and these might obstruct bijectivity of the classical
time $t_x$ defined as before in~(\ref{eq:class-time-tau}).
Indeed the derivative $t_x^\prime(\tau)=x(\tau)^2/\norm{x}^2\ge 0$, 
more properties in~(\ref{eq:class-time-tau-deriv}), might now be zero
at some times -- precisely the times of collision of the solution $x$
with the origin. By continuity of the map $x\colon[0,1]\to\R$ the zero set
$T_x:=x^{-1}(0)$ is closed, thus compact.
On the other hand, the set is discrete because being a critical point 
$x\not\equiv 0$ solves a second order delay equation, see
Corollary~\ref{cor:crit-Bb}.
We denote the finite set of \textbf{regularized collision times} by
\[
   \Tt_x:=\{\tau\in[0,1]\mid x(\tau)=0\}=\{\tau_1,\dots,\tau_N\}
\]
Thus $t_x\colon[0,1]\to[0,1]$ is still strictly increasing, hence a
bijection. We denote the inverse again by $\tau_x:={t_x}^{-1}$
and the set of \textbf{classical collision times}, this terminology will become
clear in a moment, by $T_{q_x}:=t_x(\Tt_x)$, that is
\[
   T_{q_x}:=\{t_1,\dots,t_N\},\quad t_i:=t_x(\tau_i)
\]
The derivative of $\tau_x\colon[0,1]\to[0,1]$ is still given by
\begin{equation}\label{eq:t_q^prime}
   \dot\tau_x(t)=\frac{\norm{x}^2}{x(\tau_x(t))^2}
\end{equation}
but now only at \emph{non-collision times} $t$,
that is $t\in\SS^1\setminus T_{q_x}$.

\medskip
The following proposition is a special case of a theorem
due to Barutello, Ortega, and Verzini~\cite{Barutello:2021b}.

\begin{proposition}\label{prop:reg=phys}
Given a critical point $x\in W^{1,2}_\times$ of $\Bb$, namely a
solution of~(\ref{eq:Crit-Lag}), then the
rescale-squared map $q=q_x:=\Rr(x)$ is a \textbf{physical solution}
in the sense that $q$ solves the free fall equation at all times
\[
   \ddot q(t)=-\frac{1}{q(t)^2}
   ,\qquad
   \forall t\in\SS^1\setminus T_{q_x}
\]
\emph{except at the finitely many collision times} which form the set
$T_{q_x}=\{t_1,\dots,t_N\}$.
\end{proposition}

Think of the critical points $x$ of $\Bb$
as the regularized versions, ``the regularizations'',
of the physical solutions $q_x$.
Note that when the regularization $x$ runs through the big mass sitting at
the origin the physical solution $q$ bounces back.

\begin{figure}[h]
  \centering
  \includegraphics%[width=0.9\textwidth]
%                             [height=4cm]
                             {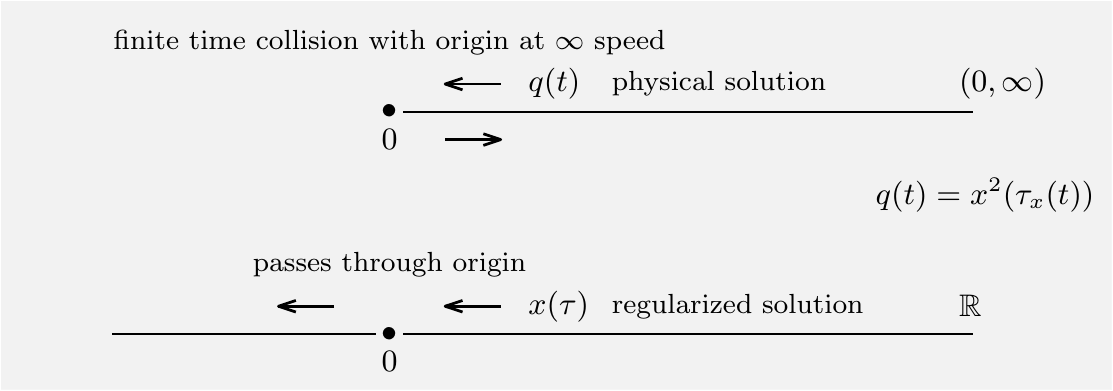}
  \caption{Physical solution $q$ and regularized solution $x$}
   \label{fig:fig-phys-traj}
\end{figure} 

\begin{proof}
At $t\in\SS^1\setminus T_{q_x}$, set $\tau(t):=\tau_x(t)$,
the derivative of $q=q_x$ is given by
\begin{equation}\label{eq:x'}
   \dot q (t)
   =2\,x(\tau(t))\, x^\prime(\tau(t)) \underbrace{\dot \tau(t)}_{(\ref{eq:t_q^prime})}
   =2\norm{x}^2\,\frac{x^\prime(\tau(t))}{x(\tau(t))}
\end{equation}
and via a change of variables from $t$ to $s:=\tau(t)$
we get using~(\ref{eq:t_q^prime}) that
\begin{equation}\label{eq:x'^2}
   \norm{\dot q}^2
   ={\color{brown} 4\norm{x}^2}\norm{x^\prime}^2
   \qquad\text{equivalently}\qquad
   \INNER{\dot q_x}{\dot q_x}=\INNER{x^\prime}{x^\prime}_x
\end{equation}
To calculate the second derivative $\ddot q$
use the critical point equation~(\ref{eq:Crit-Lag}) to get
\begin{equation}\label{eq:x''}
\begin{split}
   \ddot q(t)
   &=\frac{2\norm{x}^4}{x(\tau(t))^3}\left(x^{\prime\prime}(\tau(t))
   -\frac{x^\prime(\tau(t))^2}{x(\tau(t))}\right)\\
   &=
   \frac{1}{q(t)}\left(2\norm{x}^2\norm{x^\prime}^2-\frac{1}{\norm{x}^2}\right)
   -\frac{\dot q(t)^2}{2q(t)}\\
   &= 
   \frac{1}{q(t)}
   \left(\frac12
   \norm{\dot q}^2
   -\int_0^1 \frac{ds}{q(s)}
   -\frac12 \dot q(t)^2
   \right)
\end{split}
\end{equation}
To obtain equation three we used the identities~(\ref{eq:x'^2}) and~(\ref{eq:1/q^2}).
Let $t_-,t_+\in T_{q_x}$ be neighboring collision times in the sense that
the interval $(t_-,t_+)$ lies in the complement of the collision
time set $T_{q_x}$.
\\
Similarly to Barutello, Ortega, and Verzini~\cite[Eq.\,(3.12)]{Barutello:2021b}
we consider
\[
   \beta=\beta_q:=\frac{\ddot q}{q}
\]
as a function on the open non-collision interval $(t_-,t_+)$.
Differentiate the identity
$
   q^2\beta
   =
   \left(\frac12
   \norm{\dot q}^2
   -\int_0^1 \frac{dt}{q(t)}
   -\frac12 \dot q^2
   \right)
$
to obtain that 
$
   2q\dot q\beta+q^2\dot \beta
   =-\dot q\ddot q
   =-q\dot q\beta
$
or equivalently
$
   q^2\dot \beta=-3\beta q\dot q
$.
Hence the logarithmic derivatives satisfy
\[
   \frac{\dot \beta}{\beta}
   =-\frac{3\dot q}{q}
\]
and therefore
\[
   \beta=\frac{\mu}{q^3}
\]
for some constant $\mu\in\R$ that a priori might depend on
the interval $(t_-,t_+)$.
By definition of $\beta$ we conclude that
\[
   \ddot q
   =\frac{\mu}{q^2}
\]
on the interval $(t_-,t_+)$.
By~(\ref{eq:x''}) we get for the constant $\mu$ the expression
\begin{equation*}
\begin{split}
   \mu
   &= q(t)
   \left(\frac12
   \norm{\dot q}^2    -\int_0^1 \frac{ds}{q(s)}
   \right)
   -\frac12 q(t) \dot q(t)^2\\
   &= q(t)
   \left(\frac12
   \norm{\dot q}^2    -\int_0^1 \frac{ds}{q(s)}
   \right)
   -2\norm{x}^4 x^\prime(\tau(t))^2
\end{split}
\end{equation*}
where in the second equation we used~(\ref{eq:x'}).
Taking the limits $t\to t_\pm$ we get
\[
   \mu
   =-2\norm{x}^4 x^\prime (\tau(t_\pm))^2
   \le 0
\]
Thus the constant $\mu$ takes the same
value on the boundary of adjacent intervals.
Hence the constant $\mu$ is independent of the interval.

To see that $\mu=-1$ we multiply~(\ref{eq:x''}) by $q$
and use $q^{\prime\prime}=\mu/q^2$ to obtain
\[
   \frac{\mu}{q}
   =
   \frac12 \norm{\dot q}^2
   -\int_0^1 \frac{ds}{q(s)}
   -\frac12 \dot q(t)^2
\]
We take the mean value of this equation to get that
\[
   \mu \int_0^1 \frac{dt}{q(t)}
   =
   \frac12 \norm{\dot q}^2
   -\int_0^1 \frac{dt}{q(t)}
   -\frac12 \norm{\dot q}^2
   =
   -\int_0^1 \frac{dt}{q(t)}
\]
This shows that $\mu=-1$.
\end{proof}

%\newpage
%%%%%%%%%%%%%%%%%%%%%%%%%%%%%%%%%%%
%%%%%%%%%%%%%%%%%%%%%%%%%%%%%%%%%%%
%%%%%%% Section:  %%%%%%%%%%%%%%%%%%%%%
%%%%%%%%%%%%%%%%%%%%%%%%%%%%%%%%%%%
%%%%%%%%%%%%%%%%%%%%%%%%%%%%%%%%%%%
\section{Non-local Lagrangian mechanics}

\boldmath
%%%%%%%%%%%%%%%%%%%%%%%%%%%%%%%%%%%
%%%%%%% Subsection:  %%%%%%%%%%%%%%%%%%
%%%%%%%%%%%%%%%%%%%%%%%%%%%%%%%%%%%
\subsection{Non-local Lagrangian action $\Bb$ and $L^2_x$ inner product}
\unboldmath

In the novel approach~\cite{Barutello:2021b}
to the regularization of collisions discovered
recently by Barutello, Ortega, and Verzini
the change of time leads to a \emph{delayed} action functional $\Bb$.
Section~\ref{sec:ff} above explains this for the free fall.
The {\color{brown} delay}, that is the non-local term, is best
incorporated into the $L^2$ inner product on the loop space. More
precisely, we introduce a metric on the following Hilbert space take
away the origin, namely
\[
   W^{1,2}_\times
   :=W^{1,2}\setminus\{0\},\qquad
   W^{1,2}:=W^{1,2}(\SS^1,\R)
\]
Given a point $x\in W^{1,2}_\times$ and two tangent vectors
\[
\xi_1,\xi_2\in T_x W^{1,2}_\times=W^{1,2}
\]
we define the \textbf{\boldmath$L^2_x$ inner product} by
\begin{equation}\label{eq:q-metric}
   \inner{\xi_1}{\xi_2}_x
   :={\color{brown} 4\norm{x}^2} \inner{\xi_1}{\xi_2}
   ,\qquad
   \text{where $\inner{\xi_1}{\xi_2}:=\int_0^1 \xi_1(\tau)\xi_2(\tau)\, d\tau$}
\end{equation}
where
$
   \norm{x}:=\sqrt{\inner{x}{x}}
$
is the $L^2$ norm associated to the $L^2$ inner product.
In our case of the 1-dimen\-sional Kepler problem 
the functional then attains the form
\begin{equation*}
\begin{split}
   \Bb\colon W^{1,2}_\times:=W^{1,2}(\SS^1,\R)\setminus \{0\}\to(0,\infty)
   ,\quad
   x\mapsto \tfrac12\INNER{x^\prime}{x^\prime}_x+\frac{1}{\norm{x}^2}
\end{split}
\end{equation*}
One might interpret this functional as a non-local mechanical system
consisting of kinetic minus potential energy.

%.\newpage
\boldmath
%%%%%%%%%%%%%%%%%%%%%%%%%%%%%%%%%%%
%%%%%%% Subsection:  %%%%%%%%%%%%%%%%%%
%%%%%%%%%%%%%%%%%%%%%%%%%%%%%%%%%%%
\subsection{Critical points and  Hessian operator}
\unboldmath

Straightforward calculation provides

\begin{lemma}[Differential of $\Bb$]\label{le:diff-Bb}
The differential
$
   d\Bb\colon W^{1,2}_\times \times W^{1,2}\to\R
$
is
\begin{equation*}
\begin{split}
   d\Bb(x,\xi)
   &=4\inner{x}{\xi}\norm{x^\prime}^2+4\norm{x}^2\inner{x^\prime}{\xi^\prime}
   -2\frac{\inner{x}{\xi}}{\norm{x}^4}\\
   &=-\inner{x^{\prime\prime}}{\xi}_x
   +\underbrace{\left(\frac{\norm{x^\prime}^2}{\norm{x}^2}-\frac{1}{2\norm{x}^6}\right)}
                         _{=:\alpha<0}
   \inner{x}{\xi}_x
\end{split}
\end{equation*}
where identity two is valid whenever $x$ is of better regularity
$W^{2,2}(\SS^1,\R)\setminus \{0\}$.
\end{lemma}

Note that $\alpha=\alpha_x<0$.
Indeed $\alpha>0$ is impossible by the periodicity requirement for
solutions and $\alpha=0$ is impossible, because otherwise $x$ would
vanish identically which is excluded by assumption.

\begin{corollary}[$\Crit\,\Bb$]\label{cor:crit-Bb}
The set of critical points consists of the smooth solutions
$x\in C^\infty(\SS^1,\R)\setminus \{0\}$ of the second order delay
equation $x^{\prime\prime}=\alpha x$, see~(\ref{eq:Crit-Lag}).
\end{corollary}

By Lemma~\ref{le:diff-Bb} the
\textbf{\boldmath$L^2_x$ gradient of $\Bb$} at a loop
$x\in W^{2,2}_\times$ is given by
\begin{equation}\label{eq:grad-Bb}
   {\grad}^x\Bb(x)
   =-x^{\prime\prime}+\alpha x,\qquad
   x\in W^{2,2}_\times:= W^{2,2}(\SS^1,\R)\setminus \{0\}
\end{equation}

\boldmath
%%%%%%% Subsubsection:  %%%%%%%%%%%%%%%%
\subsubsection*{Solutions of the critical point equation}
\unboldmath

To find the solutions $x:\SS^1\to\R$ of equation~(\ref{eq:Crit-Lag})
we set $\beta:=-\alpha>0$ to get the second order ode and its solution
\[
   x^{\prime\prime}=-\beta x,\qquad
    x(\tau)=C\cos\sqrt{\beta} \tau+D\sin\sqrt{\beta} \tau
\]
where $C$ and $D$ are constants.
Thus
\[
   x^\prime=-C\sqrt{\beta}\sin\sqrt{\beta} \tau+D\sqrt{\beta}\cos\sqrt{\beta} \tau
\]
Since our solutions $x$ is periodic with period $1$
we must have
\begin{equation}\label{eq:ab-heat}
    \sqrt{\beta}=2\pi k,\quad k\in\N
\end{equation}
Thus
\begin{equation*}%\label{eq:q-heat}
    x(\tau)=C\cos 2\pi k \tau+D\sin 2\pi k \tau
\end{equation*}
and
\begin{equation*}%\label{eq:p-heat}
   x^\prime(\tau)=-2\pi k C\sin 2\pi k \tau +2\pi k D\cos 2\pi k \tau
\end{equation*}
Note that integrating the identity $1=\cos^2+\sin^2$ we get
\begin{equation}\label{eq:1/2-heat}
   1=\int_0^1 1\,d\tau
   =\int_0^1 \cos^2 \,d\tau+\int_0^1 \sin^2\,d\tau
   =2 \int_0^1 \cos^2(2\pi k \tau)\,d\tau
\end{equation}
as is well known from the theory of Fourier series.
Moreover, from the theory of Fourier series it is known
that cosine is orthogonal to sine, that is
\[
   0=\int_0^1 \cos(2\pi k \tau) \sin(2\pi k \tau)\,d\tau
\]
Therefore $\norm{x}^2=\frac{C^2+D^2}{2}$ and 
$
   \norm{x^\prime}^2
   =\frac{C^2+D^2}{2}\left(2\pi k\right)^2
$
and so we get that
\begin{equation}\label{eq:C-D-heat}
   (2\pi k)^2
   \stackrel{(\ref{eq:ab-heat})}{=}\beta^2
   \stackrel{(\ref{eq:Crit-Lag})}{=}\frac{1}{2\norm{x}^6}
   -\frac{\norm{x^\prime}^2}{\norm{x}^2}
   =\frac{4}{(C^2+D^2)^3}-(2\pi k)^2,\quad k\in\N
\end{equation}
We fix the parametrization of our solution by requiring
that at time zero the solution is maximal. Therefore
$D=0$ and we abbreviate $c_k:=C(k)>0$. Then $c_k$ is uniquely
determined by $k$ via the above equation which becomes
\begin{equation}\label{eq:c_k-0-heat}
   \frac{4}{c_k^6}
   =8\pi^2 k^2,\qquad
   \frac{1}{c_k^6}
   =2(\pi k)^2
\end{equation}
Hence
\begin{equation}\label{eq:c_k-heat}
    c_k=\frac{1}{2^\frac{1}{6}(\pi k)^\frac{1}{3}}
   \in(0,1)
   ,\qquad
    c_k^2=\frac{1}{2^\frac{1}{3}(\pi k)^\frac{2}{3}}
\end{equation}
Note that $c_k<1$. This proves

\begin{lemma}[Critical points of $\Bb$]\label{le:crit-pts-Bb}
The maps
\begin{equation}\label{eq:q-heat}
   x_k(\tau)=c_k\cos 2\pi k\tau,\qquad k\in\N
\end{equation}
are solutions of~(\ref{eq:Crit-Lag}).
Each map $x_k$ is worth a circle of solutions via time
shift $\sigma_* x_k:=x_k(\cdot+\sigma)$ where $\sigma\in\SS^1$.
All solutions of~(\ref{eq:Crit-Lag}) are given by
\[
   \Crit\,\Bb=\bigcup_{k\in\N} \{\sigma_* x_k\mid\sigma\in\SS^1\}
\]
\end{lemma}

\boldmath
%%%%%%% Subsubsection:  %%%%%%%%%%%%%%%%
\subsubsection*{The Hessian operator with respect to the $L^2_x$ inner
product}
\unboldmath

The \textbf{Hessian operator} $A_x$ of the Lagrange functional $\Bb$
is the derivative of the $L^2_x$ gradient equation $0=-x^{\prime\prime}+\alpha x$
at a critical point $x$.
Varying this equation with respect to $x$ in direction $\xi$ we obtain that
\begin{equation*}
\begin{split}
   A_x\xi
   &=-\ddot\xi+\alpha \xi
   +\left(-\frac{2\inner{x^{\prime\prime}}{\xi}}{\norm{x}^2}
   -\frac{2\norm{x^{\prime}}^2\inner{x}{\xi}}{\norm{x}^4}
   +\frac{3\inner{x}{\xi}}{\norm{x}^8}\right) x\\
   &=-\ddot\xi+\alpha \xi
   -\left(\frac{2\alpha}{\norm{x}^2}
   +\frac{2\norm{x^{\prime}}^2}{\norm{x}^4}
   -\frac{3}{\norm{x}^8}\right) \INNER{x}{\xi}  x\\
   &=-\ddot\xi+\alpha \xi
   -\frac{2}{\norm{x}^2}\left(2\alpha-\frac{1}{\norm{x}^6}\right)
   \INNER{x}{\xi}  x
\end{split}
\end{equation*}
where in step 2 we used the critical point equation 
$x^{\prime\prime} =\alpha x$ and in step 3 we replaced
$2\norm{x^{\prime}}^2$ by $2\alpha\norm{x}^2+\frac{1}{\norm{x}^4}$
according to the definition of $\alpha$.
This proves

\begin{lemma}
The Hessian operator of $\Bb$ at a critical point $x$ is given by
\begin{equation*}
\begin{split}
   A_x\colon W^{2,2}(\SS^1,\R)&\to L^2(\SS^1,\R)
   \\
   \xi&\mapsto 
   -\ddot\xi+\alpha \xi
   -\frac{2}{\norm{x}^2}\left(2\alpha-\frac{1}{\norm{x}^6}\right)
   \INNER{x}{\xi}  x
\end{split}
\end{equation*}
\end{lemma}
Recall that by~(\ref{eq:q-heat}) the critical points of $\Bb$ 
are of the form
\begin{equation}\label{eq:q_k-2}
   x_k(\tau) =c_k\cos 2\pi k\tau,\qquad k\in\N
\end{equation}
with $c_k$ given by~(\ref{eq:c_k-heat}).
Taking two $\tau$ derivatives we conclude that
\[
   x^{\prime\prime}_k=-(2\pi k)^2 x_k
\]
Since $x^{\prime\prime}_k=\alpha x_k$ we obtain
\[
   \alpha=\alpha({x_k})
   =-(2\pi k)^2
   =-\frac{2}{c_k^6}
\]
where the last equality is~(\ref{eq:c_k-0-heat}).
The formula of the Hessian operator $A_{x_k}$ involves the $L^2$ norm
of $x_k$ and, in addition, the formula of the non-local Lagrange functional $\Bb$
involves $\norm{x^{\prime}_k}^2$.
By~(\ref{eq:q-heat}) and~(\ref{eq:1/2-heat}) we obtain that
\begin{equation*}
   \Norm{x_k}^2
   =\frac{c_k^2}{2}
   =\frac{1}{2^\frac{4}{3} (\pi k)^\frac{2}{3}}
   ,\qquad
   \Norm{x^{\prime}_k}^2
   =(2\pi k)^2\frac{c_k^2}{2}
   =2^\frac{2}{3} (\pi k)^\frac{4}{3}
\end{equation*}
Thus
\[
   \Bb(x_k)
   =2 \norm{x_k}^2\norm{x^{\prime}_k}^2+\frac{1}{\norm{x_k}^2}
   =2^\frac{1}{3} 3(\pi k)^\frac{2}{3}
\]
To calculate the formula of $A_{x_k}$ we write $\xi$ as a Fourier series
\[
   \xi
   =\xi_0+\sum_{n=1}^\infty\left(
   \xi_n\cos 2\pi n\tau+\xi^n\sin 2\pi n\tau
   \right)
\]
and we use the orthogonality relation
\[
   \inner{\cos 2\pi n\cdot}{\xi}=\frac12\xi_n
\]
to calculate the product
\[
   \inner{x_k}{\xi}
   =\tfrac12 c_k\xi_k
\]
Putting everything together we obtain

\begin{lemma}[Critical values and Hessian]
The critical points of $\Bb$ are of the form
$x_k(\tau)=c_k\cos 2\pi k\tau$ for $k\in\N$,
see~(\ref{eq:q_k-2}).
At any such $x_k$ the value of~$\Bb$~is
\[
   \Bb(x_k)
   =2^\frac{1}{3} 3(\pi k)^\frac{2}{3}
\]
and the Hessian operator of $\Bb$ is
\begin{equation*}
\begin{split}
  A_{x_k}\xi
   =-\ddot \xi-(2\pi k)^2\xi
   + {\color{brown} 12(2\pi k)^2 \xi_k}\cos 2\pi k\tau
\end{split}
\end{equation*}
for every $\xi\in W^{2,2}(\SS^1,\R)$.
\end{lemma}

%%%%%%%%%%%%%%%%%%%%%%%%%%%%%%%%%%%
%%%%%%% Subsection:  %%%%%%%%%%%%%%%%%%
%%%%%%%%%%%%%%%%%%%%%%%%%%%%%%%%%%%
\subsection{Eigenvalue problem and Morse index}

Recall that $k\in\N$ is fixed since we consider the critical point $x_k$.
We are looking for solutions of the eigenvalue problem
\[
   A_{x_k}\xi=\mu\xi
\]
for $\mu=\mu(\xi;k)\in\R$ and $\xi\in W^{2,2}(\SS^1,\R)\setminus\{0\}$.
Observe that
\begin{equation*}
\begin{split}
   -\ddot\xi
   &=\sum_{n=1}^\infty (2\pi n)^2
   \left(
      \xi_n\cos 2\pi n\tau+\xi^n\sin 2\pi n\tau
   \right)
   \\
   -(2\pi k)^2\xi
   &=-(2\pi k)^2\xi_0-(2\pi k)^2\sum_{n=1}^\infty\left(
   \xi_n\cos 2\pi n\tau+\xi^n\sin 2\pi n\tau
   \right)
\end{split}
\end{equation*}
Comparing coefficients in the eigenvalue equation
$A_{x_k}\xi=\mu\xi$ we obtain
\begin{equation*}
\begin{split}
   \cos 2\pi n\tau \quad
   &
   \begin{cases}
      \mu\xi_n=4\pi^2\left(n^2-k^2\right) \xi_n
      &\text{, $\forall n\in\N_0\setminus\{k\}$}
      \\
      \mu\xi_k=
      {\color{brown} 12(2\pi k)^2 \xi_k}
      &\text{, $n= k$}
   \end{cases}
   \\
   \sin 2\pi n\tau \quad
   &
   \begin{cases}
      \mu\xi^n=4\pi^2\left(n^2-k^2\right) \xi^n
      &\text{, $\forall n\in\N$}
   \end{cases}
\end{split}
\end{equation*}

This proves

\begin{lemma}[Eigenvalues]
The eigenvalues of the Hessian $A_{x_k}$ are given by
\[
   \mu_n:=4\pi^2(n^2-k^2),\qquad n\in\N\setminus\{k\}
\]
and by
\[
   \mu_0:=-4\pi^2k^2,
   \qquad\mu_k:=0,
   \qquad\widehat\mu_k:=12(2\pi k)^2
\]
Moreover, their multiplicity (the dimension of the eigenspace) is
given by
\[
   m(\mu_n)=2,\quad n\in\N\setminus\{k\},
   \qquad m(\mu_0)=m(\mu_k)=m(\widehat\mu_k)=1
\]
\end{lemma}

Observe that the eigenvalue $\widehat \mu_k\not= \mu_n$ is different
from $\mu_n$ for every $n\in\N_0$.
Indeed, suppose by contradiction that $\widehat \mu_k=\mu_n$ for some
$n\in\N_0$, that is
\[
   48\pi^2k^2=4\pi^2(n^2-k^2)\quad\Leftrightarrow\quad
   13k^2=n^2
\]
which contradicts that $n$ is an integer.

\begin{proposition}[Morse index of $x_k$]\label{prop:Morse-calc}
The number of negative eigenvalues of the Hessian $A_{x_k}$,
called the \textbf{Morse index} of $x_k\in\Crit\,\Bb$, is odd and given by
\[
   \ind(x_k)=2k-1,\qquad k\in\N
\]
\end{proposition}

Note that by the lemma the bounded below functional $\Bb>0$ has no
critical point of index zero, in other words
\underline{the functional $\Bb>0$ has no minimum}.

%%%%%%%%%%%%%%%%%%%%%%%%%%%%%%%%%%%
%%%%%%% Section:  %%%%%%%%%%%%%%%
%%%%%%%%%%%%%%%%%%%%%%%%%%%%%%%%%%%
\section{Conley Zehnder index via spectral flow}
\label{sec:CZ}

In the local case the Conley-Zehnder index can be defined as a version
of the celebrated Maslov index~\cite{Maslov:1965a,Arnold:1967a}.
In the non-local case we do not have a flow and therefore
the interpretation as an intersection number of the linearized
flow trajectory with the Maslov cycle is not available.

In Section~\ref{sec:CZ} we recall the approach in the local case of Hofer,
Wysocki, and Zehnder~\cite{Hofer:1995b}
to the Conley-Zehnder index via winding numbers of eigenvalues
of the Hessian.
In Section~\ref{sec:nl-Ham} we will see that this theory can be
generalized to the non-local case, namely for the Hessian of the
non-local functional $\Aa_\Hh$.

\medskip
Throughout we identify the unit circle $\SS^1\subset\R^2$ with
$\R/\Z$ and $\R^2\simeq\C$ via $(x,y)\mapsto x+iy$.
Multiplication by $i$ on $\C$ is expressed by $J_0$ on $\R^2$, that is
\[
   J_0=\begin{bmatrix}0&-1\\1&0\end{bmatrix}
\]
Moreover, we think of maps $f$ with domain $\SS^1$
as $1$-periodic maps defined on $\R$, that is $f(t+1)=f(t)$, $\forall t\in\R$.

\medskip
To a continuous path\footnote{
  Because the elements of the target $L^2$ of $L_S$ are not necessarily
  continuous, a requirement to close up after time $1$ would be meaningless,
  hence paths $S:[0,1]\to \R^{2\times 2}$ are fine.
  }
of symmetric $2\times 2$ matrices
$S:[0,1]\to \R^{2\times 2}$, $t\mapsto S(t)=S(t)^T$,
we associate the operator
\begin{equation}\label{eq:L_S}
     L_S:L^2(\SS^1,\R^2)\supset W^{1,2}\to L^2(\SS^1,\R^2),\quad
     \zeta\mapsto -J_0\dot \zeta-S\zeta
\end{equation}
where $\dot\zeta:=\frac{d}{dt}\zeta$.
This operator is an unbounded self-adjoint
operator whose resolvent is compact (due to compactness
of the embedding $W^{1,2}\INTO L^2$).
Therefore the spectrum is discrete
and consists of real eigenvalues of finite multiplicity.
Given an eigenvalue $\lambda\in \spec L_S$,
then an eigenvector $\zeta$  corresponding
to $\lambda$ is a non-constantly vanishing
solution $\zeta\not\equiv 0$ of the first order ode
\[
     \dot \zeta=J_0(S+\lambda\1)\zeta
\]
for absolutely continuous $1$-periodic maps $\zeta:\R\to\R^2$.
In particular, since by definition an eigenvector is not vanishing
identically, it does not vanish anywhere, in symbols $\zeta(t)\not=0$
for every $t$.
Therefore we can associate to eigenvectors
$\zeta$ a \textbf{winding number} $w(\zeta)\in\Z$
by looking at the degree of the map
$$
     \SS^1\to\SS^1,\quad t\mapsto \frac{\zeta(t)}{\abs{\zeta(t)}} 
$$
This winding number only depends on the eigenvalue $\lambda$ and
not on the particular eigenvector for $\lambda$.
Indeed if the geometric multiplicity of $\lambda$ is $1$,
then a different eigenvector for $\lambda$ is of the form
$r\zeta$ for some nonzero real $r\not=0$.
If the geometric multiplicity of $\lambda$ is bigger than $1$,
then the space $E_\lambda$ of eigenvectors to $\lambda$ equals
a vector space of dimension at least $2$ minus the origin,
in particular, it is path connected.
The winding number has to be constant on $E_\lambda$, because
it is discrete and depends continuously on the eigenvector.
In view of these findings we write
$$
     w(\lambda)=w(\lambda;S)
$$
for the \textbf{winding number associated to an eigenvalue \boldmath$\lambda$}
of the operator $L_S$ associated to the family $S$ of symmetric matrices.

\begin{remark}[Case $S=0$]\label{rem:S=0}
Each integer $\ell\in\Z$ is realized as the winding
number of the eigenvalue $\lambda={\color{red} 2\pi}\ell$ of $L_0$ of geometric
multiplicity $2$.
\\
To see this note that the spectrum of $L_0=-J_0\frac{d}{dt}$ is
${\color{red} 2\pi}\Z$. Indeed pick an integer $\ell\in\Z$. Then the eigenvectors of
$L_0$ corresponding to ${\color{red} 2\pi}\ell$ are of the form
$t\mapsto z e^{2\pi i\ell t}$ where $z\in\C\setminus\{0\}$.
Therefore the winding number associated to $\ell$
is
\[
   w(\ell;0)=\ell
\]
The geometric multiplicity of each eigenvalue $\ell$ is $2$
(pick $z=1$ or $z=i$).
\end{remark}

\begin{remark}[General $S$]\label{rem:general-S}
Each integer $\ell\in\Z$ is realized \emph{either}
as the winding number of two different eigenvalues of $L_S$
of geometric multiplicity $1$ each
\[
     \ell=w(\lambda_1;S)=w(\lambda_2;S)
\]
\emph{or} as the winding number of a single
eigenvalue of geometric multiplicity $2$.
\\
To see this pick $\ell\in\Z$.
Consider the family $\{rS\}_{r\in[0,1]}$.
By Kato's perturbation theory we can choose 
continuous functions $\lambda_1(r)$ and $\lambda_2(r)$, $r\in[0,1]$,
such that 
\begin{enumerate}
\item[(i)]
  $\lambda_1(0)={\color{red} 2\pi}\ell$ and $\lambda_2(0)={\color{red} 2\pi}\ell$
\item[(ii)]
  $\lambda_1(r)$, $\lambda_2(r)\in\spec L_{rS}$
  and the total geometric multiplicity remains $2$, i.e.
  \[
     \dim \left(E_{\lambda_1(r)}\oplus E_{\lambda_2(r)}\right)=2
  \]
  where $E_{\lambda_i(r)}$ is the eigenspace of the eigenvalue $\lambda_i(r)$
\item[(iii)]
  $w(\lambda_1(r);rS)=\ell$ and $w(\lambda_2(r);rS)=\ell$
\end{enumerate}
The third assertion follows from the fact that
the function $r\mapsto w(\lambda_1(r);rS)$ is continuous and
takes value in the discrete set $\Z$ and, furthermore, since
$w(\ell;0)=\ell$ by Remark~\ref{rem:S=0}.
\end{remark}

Moreover, the winding number continues to be monotone in the
eigenvalue as the next lemma shows.

\begin{lemma}[Monotonicity of $w$]\label{le:monotonicity-local}
$\lambda_1<\lambda_2\in\spec L_S$ $\Rightarrow$
$w(\lambda_1)\le w(\lambda_2)$.
\end{lemma}

\begin{proof}
By contradiction suppose that there are
$\lambda_1<\lambda_2\in\spec L_S$ such that their
winding numbers satisfy
$$
     \ell_1:=w(\lambda_1;S)>w(\lambda_2;S)=:\ell_2
$$
Because eigenvalues depend continuously on the operator $L_S$
by Kato's perturbation theory, looking at the family $\{rS\}_{r\in[0,1]}$,
we can choose continuous functions $\lambda_1(r)$ and $\lambda_2(r)$
with $r\in[0,1]$ having for every $r\in[0,1]$ the following properties
\begin{enumerate}\setlength\itemsep{0ex} 
\item[(i)]
  $\lambda_1(1)=\lambda_1$ and $\lambda_2(1)=\lambda_2$
\item[(ii)]
  $\lambda_1(r)$, $\lambda_2(r)\in\spec L_{rS}$
\item[(iii)]
  $w(\lambda_1(r);rS)=\ell_1$ and $w(\lambda_2(r);rS)=\ell_2$
\end{enumerate}
From the case $S=0$ it follows that $\lambda_1(0)=w(\lambda_1(0);0)$,
which is $\ell_1$ by~(iii), and that $\lambda_2(0)=w(\lambda_2(0);0)$,
which is $\ell_2$. Thus $\lambda_1(0)=\ell_1>\ell_2=\lambda_2(0)$.
Because $\ell_1$ is different from $\ell_2$, it follows from~(iii) as
well that $\lambda_1(r)$ is different from $\lambda_2(r)$ for every $r\in[0,1]$.
Therefore by continuity we must have that
$\lambda_1(r)>\lambda_2(r)$ for every $r\in[0,1]$.
Hence by~(i) we have that $\lambda_1>\lambda_2$. Contradiction.
\end{proof}

\begin{definition}
Denote the \textbf{largest winding number among \emph{negative}
eigenvalues} of $L_S$ by
\[
     \alpha(S)
     :=\max\{w(\lambda)\mid \lambda\in (-\infty,0)\CAP \spec L_S\} \in\Z
\]
If both eigenvalues of $L_S$ (counted with multiplicity)
whose winding number is $\alpha(S)$
are $< 0$, then $S$ has \textbf{parity} $p(S):=1$.
Otherwise, define $p(S):=0$.
\end{definition}

\begin{definition}[Conley-Zehnder index of path $S$]
Following Hofer, Wysocki, and Zehnder~\cite{Hofer:1995b}
the counter-clockwise normalized Conley-Zehnder index $\CZ$ of
a continuous family of symmetric $2\times 2$ matrices
$S:[0,1]\to \R^{2\times 2}$~is
\[
     \CZ(S):=2\alpha(S)+p(S)
\]
\end{definition}

%\newpage
%%%%%%%%%%%%%%%%%%%%%%%%%%%%%%%%%%%
%%%%%%%%%%%%%%%%%%%%%%%%%%%%%%%%%%%
%%%%%%% Section:  %%%%%%%%%%%%%%%%%%%%%
%%%%%%%%%%%%%%%%%%%%%%%%%%%%%%%%%%%
%%%%%%%%%%%%%%%%%%%%%%%%%%%%%%%%%%%
\section{Non-local Hamiltonian mechanics}\label{sec:nl-Ham}

%.\newpage
%%%%%%%%%%%%%%%%%%%%%%%%%%%%%%%%%%%
%%%%%%% Subsection:  %%%%%%%%%%%%%%%%%%
%%%%%%%%%%%%%%%%%%%%%%%%%%%%%%%%%%%
\subsection{Non-local Hamiltonian equations}

Recall from~(\ref{eq:BOV-fctl}) that the functional $\Bb\colon W^{1,2}_\times\to\R$,
$x\mapsto\frac12\inner{x^\prime}{x^\prime}_x+\frac{1}{\norm{x}^2}$, extends
the classical action $\Ss_L\colon W^{1,2}_+\to\R$. Then the functional defined by
\begin{equation}\label{eq:Ll}
   \Ll\colon W^{1,2}_\times\times L^2\to\R,\quad
   (x,\xi)\mapsto\tfrac12\INNER{\xi}{\xi}_x+\frac{1}{\norm{x}^2}
\end{equation}
naturally extends $\Bb(x)=\Ll(x, x^\prime)$ in the same way
as in the classical case $\Ss_L(q)=\Ll_L(q,\dot q)$ is extended
by a corresponding functional $\Ll_L(q,v)$.

In the classical case a fiberwise strictly convex Lagrange
function $L$ on the tangent bundle determines a function $H$
on the cotangent bundle: resolve $\frak{p}:=d_{\frak{v}}L(r,\frak{v})$
for $\frak{v}$ and substitute the obtained $\frak{v}=\frak{v}(\frak{p})$
in the Legendre identity
\[
   \inner{\frak{p}}{\frak{v}}=L(r,\frak{v})+H(r,\frak{p})
\]
Returning to the non-local situation where the manifold is loop space
and $(x,\xi)$ and $(x,y)$ are pairs of loops, an analogous approach yields
\[
   y:=d_\xi\Ll(x,\xi)={\color{brown} 4\norm{x}^2} \xi,\qquad
   \xi={\color{brown}\tfrac{1}{ 4\norm{x}^2}} y,\qquad
   d_{\xi\xi}\Ll(x,\xi)=4\norm{x}^2 >0
\]
The \textbf{non-local Hamiltonian function} is then defined by
\[
   \Hh(x,y):=\inner{y}{\xi}-\Ll(x,\xi)
\]
and given by the formula
\begin{equation*}
\begin{split}
   \Hh\colon W^{1,2}_\times \times L^2&\to\R \\
   (x,y)&\mapsto
   \frac12\inner{y}{y}^x -\frac{1}{\norm{x}^2}
   =\tfrac{1}{ 4\norm{x}^2}
   \left(\tfrac{1}{2}\norm{y}^2-4\right)
\end{split}
\end{equation*}
The Hamiltonian equations of the non-local Hamiltonian $\Hh$ are the following
\begin{equation}\label{eq:ex-eqs}
  \left\{
   \begin{aligned}
     x^\prime&=\p_y\Hh  &&=\frac{1}{4\norm{x}^2}\;  y\\
     y^\prime&=-\p_x\Hh&&=\frac{x}{\norm{x}^4}\left(\frac{\norm{y}^2}{4}-2\right)\\
   \end{aligned}
   \right.
\end{equation}
for smooth solutions $x,y:\SS^1\to \R$ with $x\not\equiv 0$.

%%%%%%%%%%%%%%%%%%%%%%%%%%%%%%%%%%%
%%%%%%% Subsection:  %%%%%%%%%%%%%%%%%%
%%%%%%%%%%%%%%%%%%%%%%%%%%%%%%%%%%%
\subsection{Hamiltonian and Euler-Lagrange solutions correspond}

\begin{lemma}[Correspondence of Hamiltonian and Lagrangian solutions]
\label{le:corr-crit}\mbox{}
\begin{itemize}
\item[\rm a)]
  If $(x,y)$ solves the Hamiltonian equations~(\ref{eq:ex-eqs}), then
  $x$ solves the Lagrangian equations~(\ref{eq:Crit-Lag}).
\item[\rm b)]
  Vice versa, if $x$ solves the Lagrangian equations~(\ref{eq:Crit-Lag}),
  then $(x,y_x)$ where
  \[
     y_x:={\color{brown} 4\norm{x}^2}\,x^\prime
  \]
  solves the Hamiltonian
  equations~(\ref{eq:ex-eqs}).
\end{itemize}
\end{lemma}

\begin{proof}
a) The first equation in~(\ref{eq:ex-eqs}) leads to
$\norm{x^\prime}^2=\norm{y}^2/16\norm{x}^4$ which we resolve for $\norm{y}^2$
and then plug it into the second equation to obtain
\[
   y^\prime
   =\frac{x}{\norm{x}^4}\left(4\norm{x}^4\norm{x^\prime}^2-2\right)
\]
Now take the time $\tau$ derivative of the first equation to get indeed
\[
   x^{\prime\prime}
   =\frac{1}{4\norm{x}^2}\, y^\prime
   =\frac{1}{4\norm{x}^2}\,
   \frac{x}{\norm{x}^4}\left(4\norm{x}^4\norm{x^\prime}^2-2\right)
   =\frac{\norm{\dot x}^2}{\norm{x}^2} \, x-\frac{1}{2\norm{x}^6}\, x
   = \alpha x
\]
b) Suppose $x$ solves~(\ref{eq:Crit-Lag}) and define
$y=y_x:=4\norm{x}^2 x^\prime$, hence
$\norm{y}^2=16\norm{x}^2\norm{x^\prime}^2$.
Resolve for $x^\prime$ to obtain the first
of the Hamilton equations~(\ref{eq:ex-eqs}).
To obtain the second equation take the time $\tau$ derivative of $y$
and use the Lagrange equation for $x^{\prime\prime}$ to get
\begin{equation*}
\begin{split}
   y^\prime
   =4\norm{x}^2 x^{\prime\prime}
   =4\norm{x}^2\left(\frac{\norm{x^\prime
         }^2}{\norm{x}^2}-\frac{1}{2\norm{x}^6}\right) x
   &=\left(\frac{\norm{y}^2}{4\norm{x}^4}-\frac{2}{\norm{x}^4}\right) x\\
   &=\frac{x}{\norm{x}^4}\left(\frac{\norm{y}^2}{4}-2\right)
\end{split}
\end{equation*}
where in step three we substituted $\norm{x^\prime}^2$ according to
$\norm{y}^2=16\norm{x}^2\norm{x^\prime}^2$.
\end{proof}

\boldmath
%%%%%%%%%%%%%%%%%%%%%%%%%%%%%%%%%%%
%%%%%%% Subsection:  %%%%%%%%%%%%%%%%%%
%%%%%%%%%%%%%%%%%%%%%%%%%%%%%%%%%%%
\subsection{Non-local Hamiltonian action $\Aa_\Hh$}
\unboldmath

The \textbf{symplectic area functional} is defined by
\[
   \Aa_0\colon W^{1,2} \times L^2\to\R,\quad
   (x,y)\mapsto\int_0^1y(\tau) x^\prime (\tau)\, d\tau
\]
and the \textbf{non-local Hamiltonian action functional} by
\begin{equation*}
\begin{split}
   \Aa_\Hh:=\Aa_0-\Hh\colon 
   W^{1,2}_\times \times L^2&\to\R\\
   (x,y)&\mapsto\int_0^1y(\tau) x^\prime (\tau)\, d\tau-\Hh(x,y)
\end{split}
\end{equation*}
The derivative $d\Aa_\Hh(x,y)\colon W^{1,2}\times L^2\to\R$ is given by
\begin{equation}\label{eq:dAa}
\begin{split}
   d\Aa_\Hh(x,y)(\xi,\eta)
   &=\int_0^1\left(y\xi^\prime+\eta x^\prime\right)d\tau
   -\frac{\inner{y}{\eta}}{4\norm{x}^2}
   +\left(\tfrac{1}{4}\norm{y}^2-2\right)
   \frac{\inner{x}{\xi}}{\norm{x}^4}\\
   &=-\inner{y^\prime}{\xi}+\inner{x^\prime}{\eta}
   -\frac{\inner{y}{\eta}}{4\norm{x}^2}
   +\left(\tfrac{1}{4}\norm{y}^2-2\right)
   \frac{\inner{x}{\xi}}{\norm{x}^4}\\
   &=
   \INNER{x^\prime-\frac{y}{4\norm{x}^2}}{\eta}
   +\INNER{-y^\prime+\frac{\tfrac{1}{4}\norm{y}^2-2}{\norm{x}^4} x}{\xi}
\end{split}
\end{equation}
and this proves

\begin{lemma}\label{le:corr-crit-2}
The critical points of $\Aa_\Hh$ are precisely the $1$-periodic solutions
of the Hamiltonian equations~(\ref{eq:ex-eqs}) of $\Hh$.
\end{lemma}

In view of Lemma~\ref {le:corr-crit} and~\ref{le:corr-crit-2}
we have a one-to-one correspondence between critical points of $\Bb$
and $\Aa_\Hh$, namely $x\mapsto (x,y_x)$. Under this correspondence
the values of both functionals coincide along critical points,
see Lemma~\ref{le:Lag-domin}.

\boldmath
%%%%%%%%%%%%%%%%%%%%%%%%%%%%%%%%%%%
%%%%%%% Subsection:  %%%%%%%%%%%%%%%%%%
%%%%%%%%%%%%%%%%%%%%%%%%%%%%%%%%%%%
\subsection{Lagrangian action $\Bb$ dominates Hamiltonian one $\Aa_\Hh$}
\unboldmath

\begin{lemma}[Lagrangian domination]\label{le:Lag-domin}
There are the identities
\begin{equation*}
\begin{split}
   \Bb(x)
   &=\Aa_\Hh(x,y)
   +\frac12\Norm{2\norm{x}\, x^\prime-\frac{y}{2\norm{x}}}^2\\
   &=\Aa_\Hh(x,y)
   +
   \frac12\Norm{{\color{brown} 4\norm{x}^2}\, x^\prime-y}^2
   {\color{brown}\tfrac{1}{ 4\norm{x}^2}}
\end{split}
\end{equation*}
for every pair of loops $(x,y)\in W^{1,2}_\times \times W^{1,2}$
\end{lemma}

\begin{proof}
Just by definition of $\Aa_\Hh$ and multiplying out the inner product we obtain
\begin{equation*}
\begin{split}
   &\Aa_\Hh(x,y)+\frac12\Norm{2\norm{x}x^\prime-\frac{y}{2\norm{x}}}^2\\
   &=\inner{y}{x^\prime}-\frac{\norm{y}^2}{8\norm{x}^2}+\frac{1}{\norm{x}^2}
   +\frac12\INNER{2\norm{x}x^\prime-\frac{y}{2\norm{x}}}
   {2\norm{x}x^\prime-\frac{y}{2\norm{x}}}\\
   &=2\norm{x}\norm{x^\prime}^2+\frac{1}{\norm{x}^2}\\
   &=\Bb(x)
\end{split}
\end{equation*}
\end{proof}

\begin{corollary}[Equal values on critical points]\label{cor:crit-values-coincide}
Both functionals coincide
  \[
     \Bb(x)=\Aa_\Hh(x,y_x),\qquad y_x:={\color{brown} 4\norm{x}^2}\, x^\prime
  \]
on critical points, i.e. the solutions $x$ of~(\ref{eq:Crit-Lag}),
equivalently $(x,y_x)$ of~(\ref{eq:ex-eqs}).
\end{corollary}

In terms of the projection $\pi$ and injection $\iota$
in~(\ref{eq:proj-inj}) the corollary tells that
\[
   \Bb=\Aa_\Hh\circ\iota
   ,\qquad
   \Bb\circ\pi=\Aa_\Hh
\]
along critical points of $\Bb$, respectively of $\Aa_\Hh$.

\boldmath
%%%%%%% Subsubsection %%%%%%%%%%%%%%%%%%
\subsubsection*{The diffeomorphism $\LL$}
\unboldmath

Given the functional $\Ll$ in~(\ref{eq:Ll}), define the non-local analogue of
the diffeomorphism introduced in~\cite[p.\,1891]{Abbondandolo:2015c},
in the local context, by the formula
\begin{equation*}
\begin{split}
   \LL\colon (T L^2_\times)_0=W^{1,2}_\times\times L^2
   &\to W^{1,2}_\times\times L^2=(T^* L^2_\times)_0
   \\
   (x,\xi)&\mapsto\left(x,d_\xi\Ll(x, x^\prime+\xi)\right)=:(x,y)
\end{split}
\end{equation*}
where $(T L^2_\times)_0=(L^2_\times)_1\times (L^2)_0$ is scale calculus notation,
cf.~\cite[\S 4]{Frauenfelder:2021b}, and
\[
   y:=d_\xi\Ll(x, x^\prime+\xi)={\color{brown} 4\norm{x}^2}\,  (x^\prime+\xi)
\]
Note that since the inverse is given by
\begin{equation*}
\begin{split}
   \LL^{-1}(x,y)
   =\left(x,\frac{y}{{\color{brown} 4\norm{x}^2}}-x^\prime\right)=:(x,\xi)
\end{split}
\end{equation*}
the solutions $(x,y)$ of the Hamiltonian equations~(\ref{eq:ex-eqs})
are zeroes of $\LL^{-1}$.

As in the ode case~\cite{Abbondandolo:2015c} also in the present delay
equation situation both functionals are related through the maps
$\iota$ and $\pi$ (Lemma~\ref{le:Lag-domin}) in the form
\begin{equation*}
\begin{split}
   \Bb\circ\pi(x,y)
   =\Aa_\Hh(x,y)+\Uu^*(x,y)
   ,\qquad
   \Uu^*(x,y):=\tfrac12\Norm{\iota(x)-y}^2{\color{brown} \tfrac{1}{4\norm{x}^2}}
\end{split}
\end{equation*}
for every $(x,y)\in W^{2,2}_\times\times W^{1,2}$.
Observe that the map $\Uu^*\ge 0$ vanishes precisely along
the critical points.

With the non-negative functional $\Uu$ defined and given by
\[
   \Uu(x,\xi):=\Uu^*\circ\LL(x,\xi)
   =\tfrac12\INNER{\xi}{\xi}_x\ge 0
\]
the functionals $\Aa_\Hh$ and $\Bb$ are related by the formula
\begin{equation*}
\begin{split}
   \Aa_\Hh\circ\LL(x,\xi)
   &=\Bb(x)-\Uu(x,\xi)\\
   &\;{\color{gray} =\tfrac12\INNER{x^\prime}{x^\prime}_x+\frac{1}{\norm{x}^2}
   -\tfrac12\INNER{\xi}{\xi}_x}
\end{split}
\end{equation*}

\boldmath
%%%%%%%%%%%%%%%%%%%%%%%%%%%%%%%%%%%
%%%%%%% Subsection:  %%%%%%%%%%%%%%%%%%
%%%%%%%%%%%%%%%%%%%%%%%%%%%%%%%%%%%
\subsection{Critical points and Hessian}
\unboldmath

Defining
\[
   a:=\frac{1}{4\norm{x}^2}>0,\qquad
   b:=\left(2-\frac{\norm{y}^2}{4}\right)\frac{1}{\norm{x}^4}
   {\color{cyan} \;>0}
\]
we get
\begin{equation}\label{eq:qp}
  \left\{
   \begin{aligned}
      x^\prime &=ay\\
      y^\prime &=-bx\\
   \end{aligned}
   \right.
\end{equation}
We prove that $b\; {\color{cyan} >0}$\,:
Since our solution has to be periodic we conclude that $b$
has to be non-negative.
We claim that $b$ has to be actually positive.
Otherwise, we have the ode $x^\prime=ay$ and $y^\prime= 0$.
Thus $x^{\prime\prime}=0$, so $x$ is linear. But as the solution
must be periodic $x$ has to be constant.
This implies that $x^\prime$ and $y$ are zero,
in particular $\norm{y}^2=0$.
Thus $b\not=0$. Contradiction. This shows that $b\; {\color{cyan} >0}$.

\smallskip
We get the second order ode and its solution
\[
   x^{\prime\prime}=-abx,\qquad
    q(\tau)=c\cos\sqrt{ab} \tau+d\sin\sqrt{ab} \tau
\]
Thus
\[
   x^\prime=-c_1\sqrt{ab}\sin\sqrt{ab} \tau+c_2\sqrt{ab}\cos\sqrt{ab} \tau
\]
and
\[
   y(\tau)=\frac{x^\prime (\tau)}{a}
   =-c_1\sqrt{b/a} \sin\sqrt{ab} \tau+c_2\sqrt{b/a} \cos\sqrt{ab} \tau
\]
Since our solutions $x,y$ are periodic with period $1$
we must have
\begin{equation}\label{eq:ab}
    \sqrt{ab}=2\pi k,\quad k\in\N
\end{equation}
Thus
\begin{equation}\label{eq:q}
    x(\tau)=C\cos 2\pi k \tau+D\sin 2\pi k \tau
\end{equation}
where $C$ and $D$ are constants and
\begin{equation}\label{eq:p}
   y(\tau)=-C\frac{2\pi k}{a}\sin 2\pi k \tau +D\frac{2\pi k}{a}\cos 2\pi k \tau
\end{equation}
Note that integrating the identity $1=\cos^2+\sin^2$ we get
\begin{equation}\label{eq:1/2}
   1=\int_0^1 1\,d\tau
   =\int_0^1 \cos^2 \,d\tau+\int_0^1 \sin^2\,d\tau
   =2 \int_0^1 \cos^2(2\pi k \tau)\,d\tau
\end{equation}
as is well known from the theory of Fourier series.
Moreover, from the theory of Fourier series it is known
that cosine is orthogonal to sine, that is
\[
   0=\int_0^1 \cos(2\pi k \tau) \sin(2\pi k \tau)\,d\tau
\]
Therefore $\norm{x}^2=(C^2+D^2)/2$ and so we get for $a$ the value
\[
     a=\frac{1}{2\left(C^2+D^2\right)}
\]
Similarly we get
$
   \norm{y}^2
   =\frac{C^2+D^2}{2}\left(\frac{2\pi k}{a}\right)^2
   =2(C^2+D^2)^3 (2\pi k)^2
$
and so for $b$ we get
\[
   b=\frac{8}{\left(C^2+D^2\right)^2}-2(C^2+D^2)(2\pi k)^2
\]
Using~(\ref{eq:ab}) we get that
\begin{equation}\label{eq:C-D}
   (2\pi k)^2=ab
   =\frac{4}{\left(C^2+D^2\right)^3}-(2\pi k)^2
   ,\quad k\in\N
\end{equation}

We fix the parametrization of our solution by requiring
that at time zero the solution is maximal. Therefore
$D=0$ and we abbreviate $c_k:=C(k)>0$. Then $c_k$ is uniquely
determined by $k$ via the above equation which becomes
\begin{equation}\label{eq:c_k-0}
   \frac{4}{c_k^6}
   =8\pi^2 k^2,\qquad
   \frac{1}{c_k^6}
   =2(\pi k)^2
\end{equation}
Hence
\begin{equation}\label{eq:c_k}
    c_k=\frac{1}{2^\frac{1}{6}(\pi k)^\frac{1}{3}}
   \in(0,1)
   ,\qquad
    c_k^2=\frac{1}{2^\frac{1}{3}(\pi k)^\frac{2}{3}}
\end{equation}
Note that $c_k<1$.
Thus for each $k\in\N$ there is a solution of~(\ref{eq:qp}), namely
\begin{equation}\label{eq:qp-sol}
   \left\{
   \begin{aligned}
      x_k(\tau) &=c_k\cos 2\pi k\tau\\
      y_k(\tau) &=-2c_k^3(2\pi k)\sin 2\pi k\tau\\
   \end{aligned}
   \right.
   \qquad k\in\N
\end{equation}

%%%%%%%%%%%%%%%%%%%%%%%%%%%%%%%%%%%
%%%%%%%%%%%%%%%%%%%%%%%%%%%%%%%%%%%
\subsubsection*{Linearization}

Linearizing the Hamilton equations~(\ref{eq:ex-eqs}) at a solution
$(x,y)$ we get
\begin{equation}\label{eq:ex-lin}
  \left\{
   \begin{aligned}
     \xi^\prime&=-\frac{y}{2\norm{x}^4}\inner{x}{\xi}+\frac{\eta}{4\norm{x}^2}\\
     \eta^\prime&=\left(\frac{\xi}{\norm{x}^4}-\frac{4x}{\norm{x}^6}\INNER{x}{\xi}\right)
                             \left(\frac{\norm{y}^2}{4}-2\right)
                             +\frac12 \frac{x}{\norm{x}^4}\INNER{y}{\eta}
   \end{aligned}
   \right.
\end{equation}
for function pairs
\[
   \zeta=(\xi,\eta)\in W^{1,2}(\SS^1,\R^2)
\]
Consider the linear operator defined by
\begin{equation*}
\begin{split}
   \Ss=\Ss_{(x,y)}:W^{1,2}(\SS^1,\R^2)&\to W^{1,2}(\SS^1,\R^2)\\
   \begin{pmatrix}\xi\\\eta\end{pmatrix}
   &\mapsto
   \begin{pmatrix}
     \left(\frac{\xi}{\norm{x}^4}-\frac{4x}{\norm{x}^6}\INNER{x}{\xi}\right)
                             \left(\frac{\norm{y}^2}{4}-2\right)
      +\frac12 \frac{x}{\norm{x}^4}\INNER{y}{\eta}
     \\
     \frac{y}{2\norm{x}^4}\inner{x}{\xi}-\frac{\eta}{4\norm{x}^2}
   \end{pmatrix}
\end{split}
\end{equation*}
and the linear operator defined by
\begin{equation*}
\begin{split}
   {\color{gray} A_{(x,y)}=\;}
   L_\Ss:L^2(\SS^1,\R^2)\supset W^{1,2}&\to L^2(\SS^1,\R^2)\\
   \zeta:=\begin{pmatrix}\xi\\\eta\end{pmatrix}
   &\mapsto-J_0\zeta^\prime-I\Ss_{{\color{gray} (x,y)}}\zeta
\end{split}
\end{equation*}
where $I: W^{1,2}(\SS^1,\R^2)\INTO L^2(\SS^1,\R^2)$ is the compact
operator given by inclusion.
The kernel of the operator $L_\Ss$
is composed of the solutions to the linearized equations~(\ref{eq:ex-lin}).
If $(x,y)\in\Crit\,\Aa_\Hh$ is a critical point, then $L_\Ss$ is equal
to the Hessian operator $A_{(x,y)}$ of $\Aa_\Hh$ at $(x,y)$.

\boldmath
%%%%%%%%%%%%%%%%%%%%%%%%%%%%%%%%%%%
%%%%%%% Subsection:  %%%%%%%%%%%%%%%%%%
%%%%%%%%%%%%%%%%%%%%%%%%%%%%%%%%%%%
\subsection{Eigenvalue problem and Conley-Zehnder index}
\unboldmath

\begin{figure}%[h]
  \centering
  \includegraphics%[width=0.9\textwidth]
%                             [height=4cm]
                             {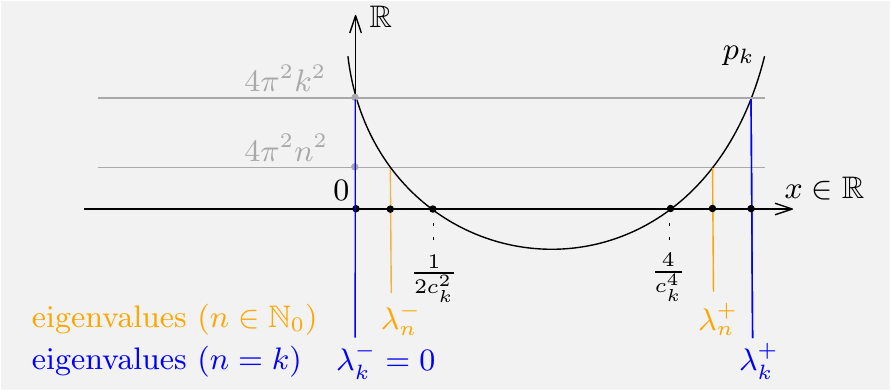}
  \caption{The parabola $p_k$ and the 
                 {\color{orange} eigenvalues $\lambda_n^\mp$}
                 for each $n\in\N_0$}
   \label{fig:fig-parabola-p_k}
\end{figure}

Fix $k\in\N$ and let $(x_k,y_k)$ be the solution~(\ref{eq:qp-sol})
of the Hamiltonian equation~(\ref{eq:ex-eqs}).
The square of the $L^2$ norm of the solution is given by
\begin{equation}\label{eq:q_k-norm}
   \Norm{x_k}^2=\frac{c_k^2}{2}
    =\frac{1}{(4\pi k)^\frac{2}{3}}
   ,\qquad
   \Norm{y_k}^2
   =2(2\pi k)^2 c_k^6
   =4
\end{equation}
Abbreviating $(x,y):=(x_k,y_k)$,
we look for reals $\lambda$ and functions
$\zeta=(\xi,\eta)$ with
\[
   L_\Ss\zeta
   :=-J_0\zeta^\prime-I\Ss\zeta
   =\lambda\zeta,\qquad
   \Ss=\Ss_{(x_k,y_k)}
\]
Apply $J_0$ to both sides of the eigenvalue problem to obtain equivalently
\begin{equation*}
\begin{split}
   \begin{pmatrix}-\lambda\eta \\ \lambda\xi\end{pmatrix}
   &=J_0\lambda \zeta
   = J_0 L_\Ss\zeta
   =\left(\p_\tau-J_0 I\Ss_{(x,y)}\right)\zeta
   \\
   &=\begin{pmatrix}
   \xi^\prime+\frac{y}{2\norm{x}^4}\inner{x}{\xi}-\frac{\eta}{4\norm{x}^2}
   \\
   \eta^\prime-\left(\frac{\xi}{\norm{x}^4}-\frac{4x}{\norm{x}^6}\INNER{x}{\xi}\right)
                            \Bigl(\underbrace{\tfrac{\norm{y}^2}{4}-2}_{=-1}\Bigr)
      -\frac12 \frac{x}{\norm{x}^4}\INNER{y}{\eta}
   \end{pmatrix}
\end{split}
\end{equation*}
Resolving for the first order terms 
and substituting $\norm{y}^2=4$
the ode becomes
\begin{equation*}
\begin{split}
   {\color{gray} \zeta^\prime=\,}
   \begin{pmatrix}\xi^\prime \\ \eta^\prime\end{pmatrix}
   &=\begin{pmatrix}
   -\lambda\eta-\frac{y}{2\norm{x}^4}\inner{x}{\xi}+\frac{\eta}{4\norm{x}^2}
   \\
   \lambda\xi-\left(\frac{\xi}{\norm{x}^4}-\frac{4x}{\norm{x}^6}\INNER{x}{\xi}\right)
      +\frac12 \frac{x}{\norm{x}^4}\INNER{y}{\eta}
   \end{pmatrix}
   {\color{gray} =J_0\lambda\zeta+J_0I\Ss\zeta}
\end{split}
\end{equation*}
Substitute first $\norm{x}^2=\frac{c_k^2}{2}$ and then $(x,y):=(x_k,y_k)$
according to~(\ref{eq:qp-sol}) to get
\begin{equation*}
\begin{split}
%   &
   \begin{pmatrix}\xi^\prime \\ \eta^\prime\end{pmatrix}
   &=\begin{pmatrix}
   -\lambda\eta-\frac{2y}{c_k^4}\inner{x}{\xi}+\frac{\eta}{2c_k^2}
     \\
   \lambda\xi-\left(\frac{4\xi}{c_k^4}-\frac{2^5x}{c_k^6}\INNER{x}{\xi}\right)
      +\frac{2x}{c_k^4}\INNER{y}{\eta}
   \end{pmatrix}
\\
   &=\begin{pmatrix}
   -\left(\lambda-\frac{1}{2c_k^2}\right)\eta
    +8\pi k\sin2\pi k\tau\inner{\cos2\pi k\,\cdot}{\xi}
     \\
   \left(\lambda-\frac{4}{c_k^4} \right)\xi
   +\frac{2^5\cos2\pi k\tau}{c_k^4}\INNER{\cos2\pi k\,\cdot}{\xi}
      -8\pi k\cos2\pi k\tau\INNER{\sin2\pi k\cdot}{\eta}
   \end{pmatrix}
\end{split}
\end{equation*}
We write the periodic absolutely continuous maps $\xi,\eta:\SS^1\to\R$
as Fourier series
\begin{equation*}
\begin{split}
   \begin{cases}
   \xi
   =\xi_0+\sum_{n=1}^\infty\left(
   \xi_n\cos 2\pi n\tau+\xi^n\sin 2\pi n\tau
   \right)
   \\
   \eta
   =\eta_0+\sum_{n=1}^\infty\left(
   \eta_n\cos 2\pi n\tau+\eta^n\sin 2\pi n\tau
   \right)
   \end{cases}
\end{split}
\end{equation*}
We set $\xi^0=\eta^0=0$. Take the derivative to get that
\begin{equation*}
\begin{split}
   \begin{cases}
   \xi^\prime
   =\sum_{n=1}^\infty\left(
   -2\pi n\cdot \xi_n\sin 2\pi n\tau+2\pi n\cdot \xi^n\cos 2\pi n\tau
   \right)
   \\
   \eta^\prime
   =\sum_{n=1}^\infty\left(
   -2\pi n\cdot \eta_n\sin 2\pi n\tau+2\pi n\cdot \eta^n\cos 2\pi n\tau
   \right)
   \end{cases}
\end{split}
\end{equation*}
By the orthogonality relation and~(\ref{eq:1/2}) we obtain
\[
   \INNER{\cos 2\pi n \,\cdot}{\xi}
   =\frac12 \xi_n,\qquad
   \INNER{\sin 2\pi n \,\cdot}{\eta}
   =\frac12 \eta^n,\qquad n\in\N
\]
Let $n\in\N_0$. Comparing coefficients we obtain from the {\bf first equations}
above
\begin{equation*}
\begin{split}
   \sin 2\pi n\tau \quad
   &
   \begin{cases}
      -2\pi n\cdot\xi_n
      =-\left(\lambda-\frac{1}{2c_k^2}\right)\eta^n&\text{, $n\not= k$}
      \\
      -2\pi k\cdot\xi_k
      =-\left(\lambda-\frac{1}{2c_k^2}\right)\eta^k
      +4\pi k\cdot\xi_k
      &\text{, $n= k$}
   \end{cases}
   \\
   \cos 2\pi n\tau \quad
   &
   \begin{cases}
      \hspace{.3cm} 2\pi n\cdot\xi^n
      =-\left(\lambda-\frac{1}{2c_k^2}\right)\eta_n&\hspace{1.6cm}\text{, $\forall n$}
   \end{cases}
\end{split}
\end{equation*}
and from the {\bf second equations}
\begin{equation*}
\begin{split}
   \cos 2\pi n\tau \quad
   &
   \begin{cases}
      \hspace{.25cm} 2\pi n\cdot\eta^n
      =\left(\lambda-\frac{4}{c_k^4}\right)\xi_n&\text{, $n\not= k$}
      \\
      \hspace{.25cm} 2\pi k\cdot\eta^k
      =\left(\lambda-\frac{4}{c_k^4}\right)\xi_k
      +\frac{2^4}{c_k^4}\xi_k-4\pi k\eta^k&\text{, $n=k$}
   \end{cases}
   \\
   \sin 2\pi n\tau \quad
   &
   \begin{cases}
      -2\pi n\cdot\eta_n
      =\left(\lambda-\frac{4}{c_k^4}\right)\xi^n&\hspace{2.82cm}\text{, $\forall n$}
   \end{cases}
\end{split}
\end{equation*}
Simplifying we get from the {\bf first equations}
\[
   \sin 2\pi n\tau\qquad
   \begin{cases}
      \text{\color{brown} a) } 2\pi n\cdot\xi_n
      =\left(\lambda-\frac{1}{2c_k^2}\right)\eta^n&\text{, $n\not= k$}
      \\
      \text{\color{cyan} b) } 6\pi k\cdot\xi_k
      =\left(\lambda-\frac{1}{2c_k^2}\right)\eta^k
      &\text{, $n= k$}
   \end{cases}
\]
\[
   \cos 2\pi n\tau \qquad
   \begin{cases}
      \text{\color{magenta} c) } 2\pi n\cdot\xi^n
      =-\left(\lambda-\frac{1}{2c_k^2}\right)\eta_n&\text{, $\forall n$}
   \end{cases}
\]
and from the {\bf second equations}
\[
   \sin 2\pi n\tau \qquad
   \begin{cases}
      \text{\color{magenta} d) } -2\pi n\cdot\eta_n
      =\left(\lambda-\frac{4}{c_k^4}\right)\xi^n&\text{, $\forall n$}
   \end{cases}
\]
\[
   \cos 2\pi n\tau \qquad
   \begin{cases}
      \text{\color{brown} e) } 2\pi n\cdot\eta^n
      =\left(\lambda-\frac{4}{c_k^4}\right)\xi_n&\text{, $n\not= k$}
      \\
      \text{\color{cyan} f) } 6\pi k\cdot\eta^k
      =\left(\lambda+\frac{12}{c_k^4}\right)\xi_k&\text{, $n=k$}
   \end{cases}
\]
\medskip\noindent
{\sc Eigenvalues.}
We obtain from equations {\color{magenta} c) and d)} that
\begin{equation}\label{eq:quadratic-1}
   \underbrace{\left(\lambda_n-\frac{1}{2c_k^2}\right)
   \left(\lambda_n-\frac{4}{c_k^4}\right)}
              _{\text{polynomial $p_k(x)$ in variable $x=\lambda_n$}}
   =4\pi^2n^2,\qquad n\in\N_0
\end{equation}
The polynomial $p_k(x)$ is illustrated by Figure~\ref{fig:fig-parabola-p_k}.
Equivalently to~(\ref{eq:quadratic-1}) we obtain the quadratic
equation for $\lambda_n$ given by
\[
    \lambda_n^2-\beta_k\lambda_n+\gamma_{k,n}=0
\]
where
\[
   \beta_k:=\frac{4}{c_k^4}+\frac{1}{2c_k^2}
   =\frac{8+c_k^2}{2c_k^4}
   ,\qquad
   \gamma_{k,n}
   :=\frac{2}{c_k^6}-4\pi^2n^2
   =4\pi^2\left(k^2-n^2\right)
\]
The solutions are
\[
     \lambda_n^-=\frac{\beta_k}{2}-\frac12\sqrt{\beta_k^2-4\gamma_{k,n}},\qquad
     \lambda_n^+=\frac{\beta_k}{2}+\frac12\sqrt{\beta_k^2-4\gamma_{k,n}}
\]
Note that $\gamma_{k,k}=0$ and therefore
\begin{equation}\label{eq:hjkhk}
   \lambda_k^-=0
\end{equation}
is zero as well and $\lambda_k^+=\beta_k$.
In case $n=0$ the quadratic equation~(\ref{eq:quadratic-1}) is already
factorized, so we read off
\[
     \lambda_0^-=\frac{1}{2c_k^2},\qquad
     \lambda_0^+=\frac{4}{c_k^4},\qquad 
   \lambda_0^-<\lambda_0^+
\]
both of which are real numbers. So the argument of the square root
is positive for $n=0$ and therefore for all $n$ (since $-\gamma_{k,n}$
is monotone increasing in $n$). Thus

\begin{lemma}[Monotonicity]\label{le:monotonicity}
The sequence $(\lambda_n^-)_{n\in\N_0}$ is
strictly monotone decreasing and
$(\lambda_n^+)_{n\in\N_0}$ is strictly monotone increasing.
\end{lemma}

\medskip\noindent
{\sc Eigenvectors.}
Recall that we had fixed $k\in \N$, in other words the solution
$(q_k,p_k)$ given by~(\ref{eq:qp-sol})
of the Hamiltonian equation~(\ref{eq:ex-eqs}).
We assume in addition\footnote{
  There appear new phenomena in case $n=0$.
  For instance, the geometric multiplicities of $\lambda_0^\mp$
  are $1$, as opposed to $2$ in case $n>0$.
  }
$n\not=0$, that is $n\in\N$.
Eigenvectors to the eigenvalues $\lambda_n^\pm$, notation  $u_n^\pm$,
can be found by setting $\eta_n:=1$, then according to
equation~{\color{magenta} c)} we define
$\xi^n_\pm:=-\frac{1}{2\pi n}\left(\lambda_n^\pm-1/2c_k^2\right)$.
The other Fourier coefficients we define to be equal $0$.
With these choices an eigenvector for $\lambda_n^\pm$
is given by the function
\[
   u_n^\pm:\SS^1\to\R^2,\quad
   \tau\mapsto
   \biggl(\underbrace{-\frac{1}{2\pi n}\left(\lambda_n^\pm
              -\frac{1}{2c_k^2}\right) }_{\xi^n_\pm}
   \sin 2\pi n\tau,\cos 2\pi n\tau\biggr)
\]
Note that the coefficient $\xi^n_+$ is strictly negative 
and $\xi^n_-$ is strictly positive since
\[
   \lambda_n^+-\frac{1}{2c_k^2}
   >\lambda_0^+-\frac{1}{2c_k^2}
   >\lambda_0^--\frac{1}{2c_k^2}
   =0
\]
and similarly
\[
   \lambda_n^--\frac{1}{2c_k^2}
   <\lambda_0^--\frac{1}{2c_k^2}
   =0
\]
Since $\xi^n_+>0$ we see that the eigenvector $u_n^+$ winds $n$ times
counter-clockwise around the origin, while $u_n^-$ winds $n$ times clockwise
around the origin since $\xi^n_-<0$.
Therefore the winding numbers equal $\pm n$, in symbols
\[
     w(u_n^\pm)=\pm n
\]
Note that in the ode (local) case Lemma~\ref{le:monotonicity}
would tell us that $w(\lambda_n^+)=n$, but in the non-local case
we cannot yet conclude independence of the choice of an eigenvector.

\begin{remark}[Case $n=0$]\mbox{}\\
{\sc Eigenvalue $\lambda_0^-=1/2c_k^2$.} 
In c) we choose $\eta_0:=1$ and
set all other Fourier coefficients zero.
Then $u_0^-=(0,1)$ is an eigenvector to the eigenvalue $\lambda_0^-$.
Since the function $u_0^-$ is constant, its winding number vanishes, in symbols
$
   w(u_0^-)=0
$.\smallskip
\\
{\sc Eigenvalue $\lambda_0^+=4/c_k^4$.}
By e) we can choose $\xi_0:=1$ and all other Fourier coefficients equal zero.
For these choices $u_0^+=(1,0)$ is an eigenvector to the eigenvalue $\lambda_0^+$.
By constancy the winding number is $0$, in symbols
$
   w(u_0^+)=0
$.
\end{remark}

\begin{remark}[Geometric multiplicity of eigenvalues
$\lambda_n^\pm$ is $\ge 2$ for $n\not=0,k$]
Instead of using  {\color{magenta} c) and d)}
one can use  {\color{brown} a) and e)}.
Setting $\eta^n:=1$ equation~a) motivates to define
$\xi_n^\pm:=-\frac{1}{2\pi n}\left(\lambda_n^+-1/2c_k^2\right)$.
With these choices a further eigenvector for $\lambda_n^\pm$
is given by the function
\[
   v_n^\pm:\SS^1\to\R^2,\quad
   \tau\mapsto
   \biggl(\underbrace{\frac{1}{2\pi
       n}\left(\lambda_n^\pm-\frac{1}{2c_k^2}\right)}_{\xi_n^\pm}
   \cos 2\pi n\tau,
   \sin 2\pi n\tau\biggr)
\]
We observe that, just as above, the winding number of the eigenvector
$v_n^+$ is $n$, and of $v_n^-$ it is $-n$, in symbols
\[
     w(v_n^\pm)=\pm n
\]
\end{remark}

%%%%%%%%%%%%%%%%%%%%%%%%%%%%%%%%%%%
%%%%%%%%%%%%%%%%%%%%%%%%%%%%%%%%%%%
\subsubsection*{Case \boldmath$n=k$ and the eigenvalues $\widehat\lambda_k^\pm$}

\begin{figure}%[h]
  \centering
  \includegraphics%[width=0.9\textwidth]
%                             [height=4cm]
                             {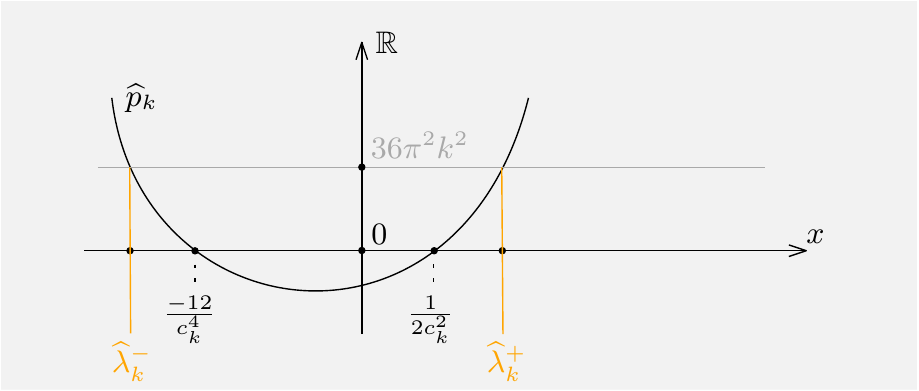}
  \caption{Parabola $\widehat p_k$ in the variable $x=\widehat\lambda_k$
                 given by~(\ref{eq:quadratic-2})}
  \label{fig:fig-p_k-hat}
\end{figure}

In the case $n=k$
we obtain from equations  {\color{cyan} b) and f)} the quadratic equation
\begin{equation}\label{eq:quadratic-2}
   \underbrace{\left(\widehat\lambda_k-\frac{1}{2c_k^2}\right)
   \left(\widehat\lambda_k+\frac{12}{c_k^4}\right)}
        _{\text{polynomial $\widehat p_k(x)$ in variable $x=\widehat\lambda_k$}}
   =36\pi^2k^2
\end{equation}
The polynomial $\widehat p_k(x)$ is illustrated by Figure~\ref{fig:fig-p_k-hat}.
{\sc Eigenvalues.}
Equivalently we obtain the quadratic equation for $\widehat\lambda_k$ given by
\[
    \widehat\lambda_k^2-B_k\widehat\lambda_k+C_k=0
\]
where
\begin{equation*}
\begin{split}
   B_k&:=\frac{1}{2c_k^2}-\frac{12}{c_k^4}
   =\beta_k-\frac{16}{c_k^4},\qquad\quad
   \beta_k:=\frac{4}{c_k^4}+\frac{1}{2c_k^2} =\frac{8+c_k^2}{2c_k^4}\\
   C_k
   &:=-\frac{6}{c_k^6}-36\pi^2k^2
   =-48\pi^2 k^2
\end{split}
\end{equation*}
with $c_k<1$ given by~(\ref{eq:c_k}).
The solutions are
\[
     \widehat\lambda_k^-=\frac{B_k}{2}-\frac12\sqrt{B_k^2-4C_k},\qquad
     \widehat\lambda_k^+=\frac{B_k}{2}+\frac12\sqrt{B_k^2-4C_k}
\]
{\sc Eigenvectors.}
Eigenvectors to the eigenvalues $\widehat\lambda_k^\pm$, notation  $u_k^\pm$,
can be found by setting $\eta^k:=1$, then equation b) motivates to define
\[
   \xi_k:=\frac{1}{6\pi k}\left(\widehat\lambda_k^\pm-1/2c_k^2\right)
\]
The other Fourier coefficients we define to be equal $0$.
With these choices an eigenvector for $\widehat\lambda_k^\pm$
is given by the function
\[
   v_k^\pm:\SS^1\to\R^2,\quad
   \tau\mapsto
   \left(\frac{1}{6\pi k}\left(\widehat\lambda_k^\pm-\frac{1}{2c_k^2}\right) \cos 2\pi k\tau,
   \sin 2\pi k\tau\right)
\]
As one sees from Figure~\ref{fig:fig-p_k-hat}
the following inequalities hold
\[
   \widehat\lambda_k^-<-\frac{12}{c_k^4}<0
    <\frac{1}{2c_k^2}<\widehat\lambda_k^+
\]
Therefore $\widehat\lambda_k^+-\frac{1}{2c_k^2}>0$ 
and $\widehat\lambda_k^--\frac{1}{2c_k^2}<0$ and hence
the winding number are
\[
     w(v_k^\pm)=\pm k
\]

\boldmath
%%%%%%%%%%%%%%%%%%%%%%%%%%%%%%%%%%%
%%%%%%% Subsection:  %%%%%%%%%%%%%%%%%%
%%%%%%%%%%%%%%%%%%%%%%%%%%%%%%%%%%%
\subsection{Disjoint families and winding numbers}
\unboldmath

\begin{figure}%[h]
  \centering
  \includegraphics%[width=0.9\textwidth]
%                             [height=4cm]
                             {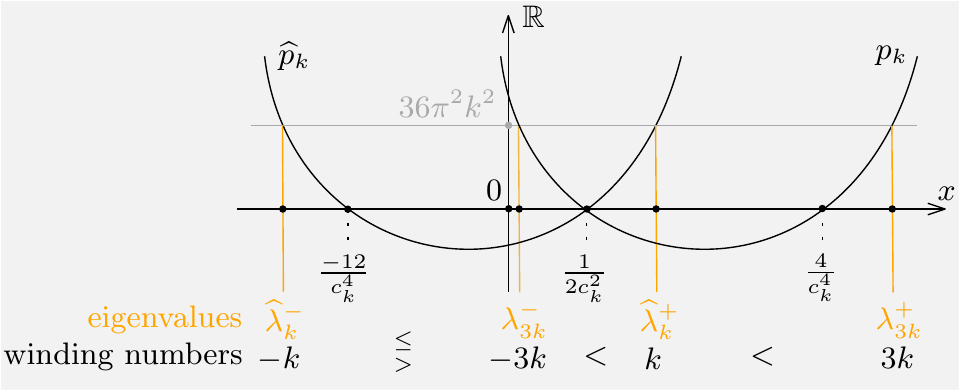}
  \caption{The parabolas $\widehat p_k$ and $p_k$, {\color{orange} eigenvalues} and
                 {\color{black} winding numbers}}
  \label{fig:fig-parabolas-dois}
\end{figure}

Consider the two quadratic polynomials $p_k$ and $\widehat p_k$ in the
variable $x=\lambda_n$
given by the left hand sides of~(\ref{eq:quadratic-1})
and~(\ref{eq:quadratic-2}), namely
\[
   p_k(x)
   :=\left(x-\frac{1}{2c_k^2}\right)\left(x-\frac{4}{c_k^4}\right)
\]
and
\[
   \widehat p_k(x)
   :=\left(x-\frac{1}{2c_k^2}\right)
   \left(x+\frac{12}{c_k^4}\right)
\]
These two polynomials have a common zero at $x=1/2c_k^2$,
they are sketched in Figure~\ref{fig:fig-parabolas-dois}.
For $n=3k$ there is equality $4\pi^2n^2=36\pi^2k^2$
and the intersection of $\widehat p_k$ and $p_k$ with the horizontal line $\{36\pi^2k^2\}$
consists of 4 points whose $x$-coordinates are the following
eigenvalues in the following order
\[
   \widehat\lambda_k^-
   <\lambda_{3k}^-
   <\widehat\lambda_k^+
   <\lambda_{3k}^+
\]

\begin{proposition}\label{prop:disj-fams}
For any $n\in\N_0$ the $\lambda_n^\mp$ are different from
$\widehat\lambda_k^-$ and from $\widehat\lambda_k^+$, in symbols
\[
   \lambda_n^\mp\not=\widehat\lambda_k^-,\qquad
   \lambda_n^\mp\not=\widehat\lambda_k^+,\qquad
   n\in\N_0
\]
\end{proposition}

\begin{proof}
Note that $\widehat\lambda_k^-<0$ is negative. On the other hand
$\lambda_n^+>0$, for $n\in\N_0$, as well as $\lambda_0^->0$ are
all positive. Therefore it suffices to show that 
$\widehat\lambda_k^-\not=\lambda_n^-$ for every $n\in\N$.

Suppose by contradiction that there are $i,j\in\{+,-\}$ such that
$\widehat\lambda_k^i=\lambda_n^j=:\lambda$
for some $n\in\N_0$. The idea is to construct two polynomials $P(z)$ and $Q(z)$
which have a common zero at $z=c_k^2$ and then use algebra to
show that there cannot be two such \smallskip polynomials.
\\
{\sc Step 0.} It is useful to consider the field extension $\Q(\pi)$ of $\Q$
which is a subfield of $\R$, that is $\Q\subset\Q(\pi)\subset\R$.
Elements of $\Q(\pi)$ have the following form. Given rational polynomials
$p,q\in\Q[x]$ with $q\not\equiv 0$ not the zero polynomial,
numbers in $\Q(\pi)$ are given by $p(\pi)/q(\pi)$.
Note that by the theorem of Lindemann~\cite{Lindemann:1882a},
see~\cite{Hilbert:1893a} for an elegant proof by Hilbert, the number
$\pi$ is transcendental\footnote{
  A real number is called {\bf transcendental}
  if it is not a zero of a polynomial with rational coefficients.
  Transcendental implies irrational. Note that $\sqrt{2}$ is
  irrational, but not transcendental (zero of $x^2-2$).
  }
and therefore $q(\pi)\not= 0$ is \smallskip non-zero.
\\
{\sc Step 1.}  The definition of the polynomial
\begin{equation}\label{eq:gjhgggj}
   Q(z):=z^3+a_0,\qquad a_0:=-(c_k^2)^3
   \stackrel{(\ref{eq:c_k})}{=}-\frac{1}{2\pi^2 k^2}\in\Q(\pi)
\end{equation}
is motivated by the goal that it has a zero at the point
$z=c_k^2$ given by~(\ref{eq:c_k}).
  Here $\Q(\pi)$ is the field extension of $\Q$ by adjoining $\pi$
  from Step 0.
  The field extension $\Q(\pi)$ is isomorphic to the field $\Q(x)$ of rational
  functions.
\\
{\sc Step 2.}  Divide the polynomial identity $\widehat p_k(\lambda)=36\pi^2 k^2$ by
$p_k(\lambda)=4\pi^2 n^2$ to get that
\[
   \frac{\lambda+\frac{12}{c_k^4}}{\lambda-\frac{4}{c_k^4}}
   =\frac{36\pi^2k^2}{4\pi^2n^2}
   =9\frac{k^2}{n^2}
\]
Now multiply by the denominator to get that
\[
   \lambda+\frac{12}{c_k^4}
   =9\frac{k^2}{n^2}\left(\lambda-\frac{4}{c_k^4}\right)
\]
Resolving for $\lambda$ yields
\[
   \lambda
   =\frac{4}{c_k^4}\cdot\frac{9\frac{k^2}{n^2}+3}{9\frac{k^2}{n^2}-1}
   =\frac{4}{c_k^4}\cdot\frac{9k^2+3n^2}{9k^2-n^2}
\]
Consequently
\begin{equation}\label{eq:gjhgj}
   \lambda\in \frac{1}{c_k^4} \Q
\end{equation}
Now evaluate $\widehat p_k$ at $\lambda$ to get
(since $\lambda:=\widehat\lambda_k^i$) that
\begin{equation}\label{eq:hu656}
\begin{split}
   0
  &=\widehat p_k(\lambda)-36 \pi^2k^2\\
  &=\left(\lambda-\frac{1}{2c_k^2}\right) 
   \left(\lambda+\frac{12}{c_k^4}\right)-36 \pi^2k^2\\
\end{split}
\end{equation}
Multiplication by ${c_k}^8$ and division by $-36\pi^2k^2$ leads to
\begin{equation}\label{eq:hu88656}
\begin{split}
   0
   &
   =\frac{\lambda^2c_k^8+\lambda c_k^4\left(12-c_k^2/2\right)-6c_k^2}
   {-36\pi^2k^2} +c_k^8\\
   &
   =-\frac{4}{9\pi^2 k^2}\left(\frac{9k^2+3n^2}{9k^2-n^2}\right)^2
   +\frac{1}{18\pi^2 k^2} \left(\frac{9k^2+3n^2}{9k^2-n^2}\right)
   c_k^2\\
   &\quad
   -\frac{4}{3\pi^2 k^2} \left(\frac{9k^2+3n^2}{9k^2-n^2}\right)
   +\frac{1}{18\pi^2 k^2}\left(\frac{9k^2+3n^2}{9k^2-n^2}\right) c_k^2
   +c_k^8\\
   &
   =P(c_k^2)
\end{split}
\end{equation}
where the polynomial $P$ is given by
\[
   P(z)=z^4+b_1 z+b_0,\qquad
   b_0, b_1\in\Q(\pi)
\]
and the coefficients $b_0$ and $b_1$ -- according
to~(\ref{eq:hu88656}) -- by
\begin{equation*}
\begin{split}
   b_0
   &=-\frac{4}{3\pi^2k^2}\frac{9k^2+3n^2}{9k^2-n^2}
   \left(\frac{3k^2+n^2}{9k^2-n^2} +1\right)
   =-\frac{48}{\pi^2}\frac{3k^2+n^2}{(9k^2-n^2)^2}
   <0
   \\
   b_1
   &=\frac{1}{9\pi^2 k^2}\left(\frac{9k^2+3n^2}{9k^2-n^2}\right) 
\end{split}
\end{equation*}

We define a linear polynomial in $z$ with coefficients in the field $\Q(\pi)$
by the formula
\[
   R(z):=P(z)-zQ(z)=(b_1-a_0) z+b_0
\]
Since $z=c_k^2$ is a zero of $P$ by~(\ref{eq:hu88656}) and of $Q$
by~(\ref{eq:gjhgggj}), it is a zero of the linear polynomial $R$.
\\
Since $b_0\not=0$ it follows that $b_1-a_0\not= 0$: otherwise
$R\equiv b_0\not= 0$ would not have a zero at all.
Since $0=R(c_k^2)=(b_1-a_0)c_k^2+b_0$ we get that
$c_k^2=-b_0/(b_1-a_0)\in\Q(\pi)$, that is $c_k^2$ is of the form
$p(\pi)/q(\pi)$ where $p,q\in\Q[z]$.

\medskip
We derive a contradiction: Evaluate~(\ref{eq:gjhgggj}) at
$z=\frac{p(\pi)}{q(\pi)}=c_k^2\in(0,1)$ to get
\[
   0=Q(c_k^2)=\frac{p(\pi)^3}{q(\pi)^3}-\frac{1}{2\pi^2 k^2}
\]
Multiply the identity by $q(\pi)^3\left(2\pi^2 k^2\right)$ to get that
\[
   0=\underbrace{p(\pi)^3\left(2\pi^2 k^2\right)}_{\deg=2\mod 3}
   -\underbrace{q(\pi)^3}_{\deg=0\mod 3}
\]
Consider the polynomial $s:=p^3r-q^3\in\Q[z],$ where $r(z):=2k^2z^2$.
\\
{\bf Claim.} $s\not\equiv 0$
\\
{\it Proof of claim.}
This follows from considering the degrees.
Note that 
\[
   \deg r=2 \mod 3,\qquad
   \deg p^3=0 \mod 3,\qquad \deg q^3=0 \mod 3
\]
Therefore
\[
   \deg p^3 r=\deg p^3+\deg r=2\mod 3,\qquad
   \deg q^3 =0\mod 3
\]
Hence $p^3r\not=q^3$ 
and consequently $s\not\equiv 0$.
\smallskip
This proves the claim.

In view of the claim we found a non-zero polynomial $s$ with rational
coefficients and the property that $s(\pi)=0$.
But this contradicts the theorem of Lindemann as explained earlier.\footnote{
  Lindemann showed that $\pi$ is transcendental:
  There is no non-zero polynomial with rational
  coefficients having $\pi$ as a zero.
  }
\end{proof}

\begin{corollary}[Well-defined winding number]
For every $n\in\N_0$ we define the winding numbers
\[
   w(\lambda_n^-):=-n,\quad
   w(\lambda_n^+):=n,\quad
   w(\widehat\lambda_k^-):=-k,\quad
   w(\widehat\lambda_k^+):=k
\]
In view of Proposition~\ref{prop:disj-fams}
these winding numbers are well-defined
and, in view of the discussion before, correspond
to the winding number of an arbitrary eigenvector
of the eigenvalue.
\end{corollary}

\begin{figure}[h]
\begin{center}
\begin{tabular}{lllllll}
\toprule
  $\spec L_\Ss$
  & $\lambda_n^-$
  & $\lambda_k^-=0$
  & $\lambda_0^-$
  & $\lambda_0^+$
  & $\lambda_k^+$
  & $\lambda_n^+$
\\
  
  & 
  & $\widehat\lambda_k^-$
  & 
  & 
  & $\widehat\lambda_k^+$
  & 
\\
\midrule
  $m$
  & $2$
  & $1$
  & $1$
  & $1$
  & $1$
  & $2$
\\
  $w$
  & $-n$
  & $-k$
  & $0$
  & $0$
  & $k$
  & $n$
\\
\bottomrule
\end{tabular}
\caption{Multiplicities $m$ and winding numbers $w$ of eigenvalues, $n\in\N\setminus\{k\}$}
\label{fig:magnetic-table}
\end{center}
\end{figure}

\begin{proposition}\label{prop:CZ-calc}
At a critical point $(x_k,y_k)$ of $\Aa_\Hh$, see~(\ref{eq:qp-sol}),
it holds that
\[
   \alpha(\Ss_{(x_k,y_k)})=w(\widehat\lambda_k^-)=-k
   ,\qquad
   p(\Ss_{(x_k,y_k)})=1
\]
and with the definitions $\CZ:=2\alpha+p$ and $\CZcan:=-\CZ$ we obtain
\begin{equation*}
   \CZ(x_k,y_k)=-2k+1
   ,\qquad
   \CZcan(x_k,y_k)=2k-1
\end{equation*}
\end{proposition}

\begin{proof}
With $\Ss:=\Ss_{(x_k,y_k)}$ we recall the definition of $\alpha$, namely
\[
     \alpha(\Ss)
     :=\max\{w(\lambda)\mid \lambda\in (-\infty,0)\CAP \spec L_\Ss\} \in\Z
\]
Observe that $\widehat\lambda_k^-<0$, see Figure~\ref{fig:fig-p_k-hat}, 
and $w(\widehat\lambda_k^-)=-k$,  see Figure~\ref{fig:magnetic-table}.
Therefore $\alpha(\Ss)\ge -k$.
To show the reverse inequality $\alpha(\Ss)\le -k$ we need to check
non-negativity of all eigenvalues with winding number $>-k$.
By Figure~\ref{fig:magnetic-table} the eigenvalues of winding number
$>-k$ are of three types:

\begin{itemize}\setlength\itemsep{0ex} 
\item[(i)]
  $\lambda_n^-$ for $n<k$:
  In this case $\lambda_n^->\lambda_k^-=0$ for every $n<k$
  by monotonicity, see Lemma~\ref{le:monotonicity}, and~(\ref{eq:hjkhk}).
\item[(ii)]
  $\lambda_n^+$ for all $n\in\N_0$: In this case
  $\lambda_n^+>0$ for every $n\in\N_0$ by
  Figure~\ref{fig:fig-parabola-p_k}
\item[(iii)]
  $\widehat\lambda_k^+$: In this case
  $\widehat\lambda_k^+>0$ by Figure~\ref{fig:fig-p_k-hat}.
\end{itemize}
This proves that $\alpha(\Ss)=-k$.

Because there exists a non-negative eigenvalue, namely $\lambda_k^-=0$, with
the same winding number $-k$ as the negative eigenvalue
$\widehat\lambda_k^-<0$ that realizes the maximal winding number
$\alpha(\Ss_{(x_k,y_k)})$ among negative eigenvalues, we get that
\[
   p(\Ss_{(x_k,y_k)})=1
\]
\end{proof}

%%%%%%%%%%%%%%%%%%%%%%%%%
%%%%%%%%% REFERENCES %%%%%%
%%%%%%%%%%%%%%%%%%%%%%%%
\bibliographystyle{alpha}

\addcontentsline{toc}{section}{References}
\bibliography{$HOME/Dropbox/0-Libraries+app-data/Bibdesk-BibFiles/library_math,$HOME/Dropbox/0-Libraries+app-data/Bibdesk-BibFiles/library_math_2020,$HOME/Dropbox/0-Libraries+app-data/Bibdesk-BibFiles/library_physics}{}
%$

\end{document}